\documentclass[11pt]{amsart}

\usepackage{amsmath}
\usepackage{mathrsfs}
\usepackage[all]{xy}
\usepackage{amssymb}
\usepackage{comment}
\usepackage{color}
\usepackage[colorlinks,citecolor=blue,urlcolor=black,linkcolor=black]{hyperref}

\newtheorem{theorem}{Theorem}[section]
\newtheorem{claim}[theorem]{Claim}

\newtheorem{lemma}[theorem]{Lemma}

\newtheorem{corollary}[theorem]{Corollary}

\theoremstyle{definition}
\newtheorem{definition}[theorem]{Definition}

\newtheorem{question}[theorem]{Question}

\theoremstyle{remark}

\newcount\skewfactor
\def\mathunderaccent#1#2 {\let\theaccent#1\skewfactor#2
\mathpalette\putaccentunder}
\def\putaccentunder#1#2{\oalign{$#1#2$\crcr\hidewidth
\vbox to.2ex{\hbox{$#1\skew\skewfactor\theaccent{}$}\vss}\hidewidth}}
\def\name{\mathunderaccent\tilde-3 }

\def\smallbox#1{\leavevmode\thinspace\hbox{\vrule\vtop{\vbox
   {\hrule\kern1pt\hbox{\vphantom{\tt/}\thinspace{\tt#1}\thinspace}}
   \kern1pt\hrule}\vrule}\thinspace}

\newcommand{\stick}{{\ensuremath \mspace{2mu}\mid\mspace{-12mu} {\raise0.6em\hbox{$\bullet$}}}}

\newcommand{\cf}{{\rm cf}}

\DeclareMathOperator{\crit}{crit}

\DeclareMathOperator{\image}{''}
\DeclareMathOperator{\Add}{Add}

\DeclareMathOperator{\dom}{dom}
\DeclareMathOperator{\range}{range}
\newcommand{\GCH}{\mathrm{GCH}}

\def\qedref#1{$\qed_{\reforiginal{#1}}$}

\title{On a problem of Erd\H{o}s and Hajnal}
\author{Shimon Garti}
\address{Einstein Institute of Mathematics,
 The Hebrew University of Jerusalem,
 Jerusalem 91904, Israel}
\email{shimon.garty@mail.huji.ac.il}

\author{Yair Hayut}
\address{Einstein Institute of Mathematics,
 The Hebrew University of Jerusalem,
 Jerusalem 91904, Israel}
\email{yair.hayut@mail.huji.ac.il}

\author{Saharon Shelah}
\address{Institute of Mathematics
 The Hebrew University of Jerusalem,
 Jerusalem 91904, Israel
 and  Department of Mathematics
 Rutgers University
 New Brunswick, NJ 08854, USA}
\email{shelah@math.huji.ac.il}
\urladdr{http://www.math.rutgers.edu/\char`\~shelah}
\thanks{The research was supported by Israel Science Foundation Grant no.\ 1838/19 and Grant no.\ 2320/23. The second author's research was supported by the Israel Science Foundation Grant no.\ 1963/21. This is publication 1249 (the 204th prime number) of the third author, publication 47 (the 15th prime number) of the first author and publication 37 (the 12th prime number) of the second author}


\begin{document}
\let\labeloriginal\label
\let\reforiginal\ref

\begin{abstract}
We address a question of Erd\H{o}s and Hajnal about the ordinary partition relation $\aleph_{\omega+1}\nrightarrow(\aleph_{\omega+1},(3)_{\aleph_0})^2$.
For $\theta=\cf(\lambda)<\lambda$, assuming $2^\lambda=\lambda^+$ they proved the negative relation $\lambda^+\nrightarrow(\lambda^+,(3)_\theta)^2$ and asked whether the (local instance of) \textsf{GCH} is indispensable.
We show that this negative relation is consistent with $\lambda$ being a strong limit and $2^{\lambda}>\lambda^+$.
The result can be pushed down to $\aleph_{\omega}$.
\end{abstract}

\maketitle

\newpage

\section{Introduction}

Let $G=(V,E)$ be a graph of size $\lambda$.
One may wonder whether there must be a monochromatic triangle under any edge coloring $c:E\rightarrow\theta$.
The answer is trivially no, since the graph can be a set of isolated vertices with no edges at all, or a triangle-free graph.
Thus in order to make the above question interesting, one should assume that there are many edges (and, in particular, many triangles) in the graph.
One possible way to do it uses the following concept.
A set of vertices $W\subseteq V$ is called \emph{independent} if $[W]^2\cap E=\varnothing$.
If $G$ is of size $\lambda$ and there are no independent subsets of size $\lambda$ in $G$, then there are many edges in the graph and the question makes more sense.

The above discussion can be formulated in the language of partition calculus, without mentioning graphs at all.
The ordinary partition relation $\lambda\rightarrow(\kappa,(3)_\theta)^2$ says that for every coloring $c:[\lambda]^2\rightarrow\theta$ there is either $A\in[\lambda]^\kappa$ so that $c''[A]^2=\{0\}$, or $B\in[\lambda]^3$ and $\gamma\in(0,\theta)$ so that $c''[B]^2=\{\gamma\}$.
A particularly interesting case is when $\kappa=\lambda$.
In terms of graph theory, one can interpret the coloring as assigning zero to pairs of vertices with no edge, and some color $\gamma\in(0,\theta)$ to edges of a given graph.
The positive relation $\lambda\rightarrow(\lambda,(3)_\theta)^2$ means that either there is an independent set of size $\lambda$, or a monochromatic triangle.

Erd\H{o}s, Hajnal and Rado investigated this relation in \cite{MR0202613}.
They established several results, and focused in particular on graphs whose size is a successor of a singular cardinal.
A good account appears in the monograph \cite{MR795592}, in which the following is phrased and proved:

\begin{theorem}
  \label{thmehmr20} Assume that $\lambda$ is a singular cardinal and $2^\lambda=\lambda^+$.
  Then $\lambda^+\nrightarrow(\lambda^+,(3)_{\cf(\lambda)})^2$.
\end{theorem}

Actually, they proved something a bit stronger, see \cite[Theorem 20.2]{MR795592}.
A natural question is whether the assumption $2^\lambda=\lambda^+$ is removable.
Let us indicate that if one forces $2^{\cf(\lambda)}\geq\lambda^+$ then a negative result obtains, namely $\lambda^+\nrightarrow(\lambda^+,(3)_{\cf(\lambda)})^2$, as mentioned in \cite{MR795592}.
Thus we shall assume from now on that $2^{\cf(\lambda)}<\lambda$, and in fact we shall force the negative relation while $\lambda$ is a strong limit singular cardinal.
The first case, in this context, is $\lambda=\aleph_\omega$.
In a collection of unsolved problems \cite[Problem 5]{MR0280381}, the pertinent question appears as follows:

\begin{question}
  \label{ehopen} Can one prove without assuming \textsf{GCH} that the relation $\aleph_{\omega+1}\nrightarrow(\aleph_{\omega+1},(3)_{\aleph_0})^2$ holds?
\end{question}

It appeared, again, in \cite[Problem 20.1]{MR795592}.\footnote{In the monograph \cite{MR795592}, the domain of the coloring is $\aleph_\omega^{\aleph_0}.$ Under the assumption $2^{\aleph_\omega}=\aleph_{\omega+1}$, these two entities coincide, i.e. $\aleph_\omega^{\aleph_0}=\aleph_{\omega+1}$.}
Despite the fact that powerful methods for dealing with successors of singular cardinals are available today, the problem is still open.
A recent survey by Komj\'ath, \cite{komj}, describes the progress in every problem from the list of Erd\H{o}s and Hajnal, and according to this survey no progress has been made with regard to this problem.
In this paper we intend to show that $\lambda^+\nrightarrow(\lambda^+,(3)_{\cf(\lambda)})^2$ is consistent with $2^\lambda>\lambda^+$, where $\lambda$ is singular and strong limit.
The current status of the question is, therefore, as follows.
The negative relation is consistent with the failure of \textsf{GCH}, but we do not know whether it holds in \textsf{ZFC}.

Our strategy is to replace the hypothesis $2^\lambda=\lambda^+$ by pcf arguments.
More specifically, we obtain unbounded sequences of regular cardinals below $\lambda$ that satisfy some relevant negative arrow relation, and we make sure that the true cofinality of these sequences is $\lambda^+$.
These assumptions enable us to lift the negative relations below $\lambda$ to $\lambda^+$, even if $2^\lambda>\lambda^+$.

This approach is rendered here twice.
Firstly, we get the negative arrow relation $\lambda_i^+\nrightarrow(\lambda_i^+,(3)_{\cf(\lambda_i)})^2$ where each $\lambda_i$ is a successor of a singular cardinal, simply by assuming \textsf{GCH} below $\lambda$, and then we lift this relation to $\lambda^+$.
The drawback here is that $\lambda$ must be a limit of singular cardinals with the same cofinality.
Namely, the cofinality of $\lambda$ should be the cofinality of every element in the chosen sequence of singular cardinals which tend to $\lambda$.
Thus $\aleph_\omega$ cannot be handled in this way.
Secondly, we modify the negative arrow relation over each $\lambda_i$ in such a way that it applies to regular cardinals even though these cardinals are not successors of singular cardinals.
This strategy applies to small cardinals like the $\aleph_n$'s, and consequently the negative arrow relation can be forced over $\aleph_\omega$.

There is another possible approach towards the same negative arrow relation.
This approach is based on prediction principles from the tiltan family.
Suppose that diamond holds at $\lambda^+$, in which case $2^\lambda=\lambda^+$.
In many cases, diamond can be replaced by weaker prediction principles like tiltan or stick.
The version employed in this paper is $\stick(\lambda)$, and it says that there exists a sequence $(A_\alpha\mid\alpha\in\lambda^+)$ of elements of $[\lambda^+]^\lambda$ that guesses every $A\in[\lambda^+]^{\lambda^+}$ in the sense that $A_\alpha\subseteq{A}$.
Usually, these principles at $\lambda^+$ are consistent with $2^\lambda>\lambda^+$.
We shall prove that $\stick(\lambda)$ yields the negative arrow relation $\lambda^+\nrightarrow(\lambda^+,(3)_{\cf(\lambda)})^2$.
However, we do not know whether $\stick(\lambda)$ is consistent with $2^\lambda>\lambda^+$ when $\lambda$ is a strong limit singular cardinal.
Recall that in $\stick(\lambda)$ the size of the guessing sets is $\lambda$, and this is the main problematic issue in our context.
This point will be discussed in the relevant section.

The rest of the paper contains three additional sections, and an appendix.
In the first section we solve the general problem using simple arguments of pcf theory.
As mentioned above, this solution does not cover the case of $\aleph_\omega$.
In the second section we modify the pcf arguments and show how to push down the result to $\aleph_\omega$.
Finally, in the third section we show that a special version of stick at $\lambda^+$ yields the desired negative arrow relation.
The appendix is devoted to a formal proof of the strong Prikry property with respect to the Extender-based Prikry forcing with collapses, the one that we use in our proofs.

Our notation is coherent with \cite{MR795592}.
We shall use the Jerusalem forcing notation, namely we force upwards.
A function $f:E\rightarrow\mathcal{P}(E)$ is a \emph{set mapping} if $x\notin f(x)$ whenever $x\in{E}$.
A subset $X\subseteq{E}$ is free for $f$ iff $f(y)\cap{X}=\varnothing$ whenever $y\in{X}$.
If $\kappa=\cf(\kappa)<\lambda$ then we let $S^\lambda_\kappa=\{\delta\in\lambda\mid\cf(\delta)=\kappa\}$.
If $\aleph_0<\cf(\lambda)$ then $S^\lambda_\kappa$ is a stationary subset of $\lambda$.
A cardinal $\kappa$ is $\gamma$-supercompact if there exists an elementary embedding $\jmath:V\rightarrow{M}$ such that $\kappa=crit(\jmath)$ and ${}^\gamma{M}\subseteq{M}$.
One says that $\kappa$ is supercompact if $\kappa$ is $\gamma$-supercompact for every ordinal $\gamma$.

Supercompactness can be characterized by the following (global) property.
Suppose that $\kappa\leq\lambda$.
An ultrafilter $\mathscr{U}$ over $\mathcal{P}_\kappa\lambda$ is \emph{fine} if $\{x\in\mathcal{P}_\kappa\lambda\mid\alpha\in{x}\}\in\mathscr{U}$ whenever $\alpha\in\lambda$.
One says that $\mathscr{U}$ is \emph{normal} if every $\in$-regressive function on a set $A$ from $\mathscr{U}$ is constant on a set $B\subseteq{A}, B\in\mathscr{U}$.
Using these definitions, $\kappa$ is $\lambda$-supercompact if there exists a fine and normal measure over $\mathcal{P}_\kappa\lambda$, and $\kappa$ is supercompact if $\kappa$ is $\lambda$-supercompact for every $\lambda\geq\kappa$.
We shall use the idea of indestructibility (at supercompact cardinals) as appeared in the seminal work of Laver, \cite{MR0472529}.
It is shown there that a supercompact cardinal $\kappa$ can be made indestructible under any generic extension by $\kappa$-directed-closed forcing notions.
In our context, we need a modification of Laver's proof, applicable to strategically-closed forcing notions.
This version will be phrased and proved below.

For basic background concerning Prikry type forcings we refer to \cite{MR2768695}, and to the papers of Magidor \cite{MR491183} and \cite{MR491184} in which the basic method of Prikry forcing with interleaved collapses was introduced.
We also refer to \cite{yairbd} in this context.
For background in pcf theory we suggest \cite{MR2768693}, \cite{MR1086455} and \cite{MR1318912}.

\newpage

\section{A negative relation from simple pcf arguments}

In this section we suggest our first approach for proving the negative relation $\lambda^+\nrightarrow(\lambda^+,(3)_{\cf(\lambda)})^2$ where $\lambda$ is a strong limit singular cardinal and $2^\lambda>\lambda^+$.
The idea is to assume the negative relation at an unbounded sequence of cardinals below $\lambda$ (by assuming \textsf{GCH} at these cardinals) and to obtain the negative relation at $\lambda^+$ by means of pcf theory.

Let $\mathcal{J}$ be an ideal over $\theta$, and let $(\lambda_i\mid i\in\theta)$ be a sequence of regular cardinals.
Let $(f_\alpha\mid\alpha\in\mu)$ be a sequence of functions in $\prod_{i\in\theta}\lambda_i$.
One says that this sequence is $\mathcal{J}$-increasing if $\{i\in\theta\mid f_\alpha(i)\geq f_\beta(i)\}\in\mathcal{J}$ whenever $\alpha<\beta<\mu$.
The sequence is cofinal if for every $g\in\prod_{i\in\theta}\lambda_i$ there exists $\alpha\in\mu$ so that $g\leq_{\mathcal{J}}f_\alpha$.
Occasionally, $\mathcal{J}$ will be the ideal of bounded subsets of $\theta$, denoted $J^{bd}_\theta$.
Recall that ${\rm tcf}(\prod_{i\in\theta}\lambda_i,J)=\kappa$ iff there is a $J$-cofinal and increasing sequence in the product $\prod_{i\in\theta}\lambda_i$, and $\kappa$ is the minimal length of such a sequence.
We commence with the combinatorial theorem, followed by a description of the ways to force the assumptions in this theorem.

\begin{theorem}
  \label{thmnegfrompcf} Assume that:
  \begin{enumerate}
    \item [$(a)$] $\mu>\cf(\mu)=\theta$.
    \item [$(b)$] $\mu$ is a strong limit cardinal.\footnote{As indicated by the referee, the strong limitude of $\mu$ follows from $(d)$ and $(f)$ below.}
    \item [$(c)$] $2^\mu>\mu^+$.
    \item [$(d)$] $(\mu_i\mid i\in\theta)$ is increasing and $\mu=\bigcup_{i\in\theta}\mu_i$.
    \item [$(e)$] $\mu_i>\cf(\mu_i)=\theta$ for every $i\in\theta$.
    \item [$(f)$] $\mu_i$ is a strong limit cardinal for every $i\in\theta$.
    \item [$(g)$] $2^{\mu_i}=\mu_i^+$ for every $i\in\theta$.
    \item [$(h)$] ${\rm tcf}(\prod_{i\in\theta}\mu_i^+,J^{\rm bd}_\theta)=\mu^+$.
  \end{enumerate}
  Then $\mu^+\nrightarrow(\mu^+,(3)_{\cf(\mu)})^2$.
\end{theorem}

\par\noindent\emph{Proof}. \newline
For every $i\in\theta$ let $c_i:[\mu_i^+]^2\rightarrow\theta$ be a witness to the negative relation $\mu_i^+\nrightarrow(\mu_i^+,(3)_{\theta})^2$. This negative relation follows from assumption $(g)$.
Our goal is to define a coloring $c:[\mu^+]^2\rightarrow\theta$ by combining the $c_i$s together in such a way that the corresponding negative relation at $\mu^+$ will follow.

We need two mathematical objects to define our coloring.
The first is a scale $(f_\alpha\mid\alpha\in\mu^+)$ in the product $(\prod_{i\in\theta}\mu_i^+,J^{\rm bd}_\theta)$.
The second is a system of functions $(h_i\mid i\in\theta)$ where $h_i\in{}^\theta\theta$ is injective and $h_i(0)=0$ for each $i\in\theta$.
Also, if $i<j<\theta$ then $rang(h_i)\cap rang(h_j)=\{0\}$.
Suppose that $\alpha<\beta<\mu^+$.
Let $i(\alpha,\beta)$ be the minimal $j\in\theta$ so that $f_\alpha(j)\neq f_\beta(j)$.
Such an ordinal always exists since $f_\alpha<_{J^{\rm bd}_\theta}f_\beta$.
We define the coloring $c:[\mu^+]^2\rightarrow\theta$ as follows.
Given $\alpha<\beta<\mu^+$ let $i=i(\alpha,\beta)$ and let $c(\alpha,\beta)=h_{i}(c_{i}(\{f_\alpha(i),f_\beta(i)\}))$.\footnote{The fact that $i=i(\alpha,\beta)$ implies that $\{f_\alpha(i),f_\beta(i)\}$ is a pair of ordinals.}
Let us show that $c$ exemplifies the negative relation $\mu^+\nrightarrow(\mu^+,(3)_{\cf(\mu)})^2$.

\begin{enumerate}
  \item [$(\aleph)$] Assume that $A\in[\mu^+]^{\mu^+}$.
  For every $i\in\theta$ let $A_i=\{f_\alpha(i)\mid\alpha\in{A}\}$.
  Set $X=\{i\in\theta\mid\mu_i^+=\bigcup A_i\}$, and notice that $X=\theta\ \text{mod}\ J^{\rm bd}_\theta$.
  Fix $i\in{X}$.
  For unboundedly many $\varepsilon\in{A_i}$ we choose $\alpha_\varepsilon\in{A}$ such that if $\varepsilon<\zeta$ then $\alpha_\varepsilon<\alpha_\zeta$ and $i(\alpha_\varepsilon,\alpha_\zeta)=i$.
  The choice is possible since $A_i$ is unbounded in $\mu_i^+$.
  Since $|\prod_{j\in{i}}\mu_j^+|<\mu_i^+$ (this inequality follows from $(e)$ and $(f)$), there are a set $B_i\in[\mu_i^+]^{\mu_i^+}$ and a fixed element $g\in\prod_{j\in{i}}\mu_j^+$ such that if $\varepsilon<\zeta$ are taken from $B_i$ then $\alpha_\varepsilon<\alpha_\zeta$ and $f_{\alpha_\varepsilon}\upharpoonright{i}=g$ for every $\varepsilon\in B_i$.
  Since $c_i$ witnesses the negative relation $\mu_i^+\nrightarrow(\mu_i^+,(3)_{\theta})^2$, one can choose $\varepsilon,\zeta\in{B_i}$ such that $\varepsilon<\zeta$ and $c_i(\varepsilon,\zeta)\neq{0}$.
  But then $c(\alpha_{\varepsilon},\alpha_{\zeta})\neq{0}$, so the proof of the first case is accomplished.
  \item [$(\beth)$] Assume that $\alpha<\beta<\gamma<\mu^+$.
  If $i(\alpha,\beta)\neq i(\alpha,\gamma)$ or $i(\alpha,\beta)\neq i(\beta,\gamma)$ or $i(\alpha,\gamma)\neq i(\beta,\gamma)$ then $\{\alpha,\beta,\gamma\}$ cannot be $c$-monochromatic with any color $\xi>0$ since for $i\neq{j}$ one has $rang(h_i)\cap rang(h_j)=\{0\}$ and by the definition of $c$.
  If $i(\alpha,\beta)=i(\alpha,\gamma)=i(\beta,\gamma)=i$ then $c\upharpoonright[\{\alpha,\beta,\gamma\}]^2=\{\xi\}$ with $\xi>0$ implies $c_i\upharpoonright[\{f_\alpha(i),f_\beta(i),f_\gamma(i)\}]^2=\{\xi\}$, since $h_i$ is injective.
  But this is impossible by the choice of $c_i$, so we are done.
\end{enumerate}

\hfill \qedref{thmnegfrompcf}

A corollary to the above theorem gives a partial answer to the question of Erd\H{o}s and Hajnal, by proving the consistency of the negative relation even if $2^\mu>\mu^+$.
Within the proof of this corollary we force with $\mathbb{Q}_{\bar{\mu}}$ from \cite[Definition 2.3]{MR2987137}.
For being self-contained, we unfold the definition of this forcing notion.

Let $\mu$ be supercompact, and let $\bar{\mu}=(\mu_i\mid i\in\mu)$ be an increasing sequence of regular cardinals so that $2^{|i|}<\mu_i$ for every $i\in\mu$.
A condition $p\in\mathbb{Q}_{\bar{\mu}}$ is a pair $(\eta,f)=(\eta^p,f^p)$ such that $\ell{g}(\eta)\in\mu$ and $\eta\in\prod_{i\in\ell{g}(\eta)}\mu_i$.
We refer to $\eta$ as the stem of $p$.
Also, $f\in\prod_{i\in\mu}\mu_i$ and $\eta\triangleleft{f}$.\footnote{The relation $\eta\triangleleft{f}$ means that $\eta$ is an initial segment of $f$.}
If $p,q\in\mathbb{Q}_{\bar{\mu}}$ then $p\leq{q}$ iff $\eta^p\trianglelefteq\eta^q$ and $f^p(j)\leq f^q(j)$ for every $j\in\mu$.

Intuitively, $\mathbb{Q}_{\bar{\mu}}$ adds a function $h\in\prod_{i\in\mu}\mu_i$ which dominates every old function in this product.
If $2^\mu=\mu^+$ in the ground model then $\mathbb{Q}_{\bar{\mu}}$ is $\mu^+$-cc.
Also, $\mathbb{Q}_{\bar{\mu}}$ is $(<\mu)$-strategically-closed.
Hence one can iterate $\mathbb{Q}_{\bar{\mu}}$ and preserve cardinals.
Moreover, the supercompactness of $\mu$ is also preserved as shown in \cite{MR2987137}.
If $\Upsilon=\cf(\Upsilon)>\mu$ is the length of the iteration then the generic functions added at each step form a scale.
Moreover, upon singularizing $\mu$ either by Prikry or by Magidor forcing one preserves the properties of this scale, thus forcing ${\rm tcf}(\prod_{i\in\cf(\mu)}\mu_{\rho_i},J^{\rm bd}_{\cf(\mu)})=\Upsilon$ in the generic extension, where $(\rho_i\mid i\in\cf(\mu))$ is the Prikry/Magidor generic sequence.
The true cofinality will be preserved after singularizing $\mu$, see \cite[Lemma 3.1 and Remark 3.2]{MR2987137}.

\begin{corollary}
  \label{coranswertoeh5} Assuming the existence of a supercompact cardinal in the ground model, one can force $\mu^+\nrightarrow(\mu^+,(3)_{\cf(\mu)})^2$ with $2^\mu>\mu^+$ where $\mu$ is a strong limit cardinal.
\end{corollary}

\par\noindent\emph{Proof}. \newline
Our goal is to force the assumptions of Theorem \ref{thmnegfrompcf}.
Let $\mu$ be a supercompact cardinal and fix a regular cardinal $\aleph_0\leq\theta\in\mu$.
We may assume that the supercompactness of $\mu$ is indestructible under $\mathbb{Q}_{\bar{\mu}}$, and \textsf{GCH} holds above $\mu$.
Let $(\mu_i\mid i\in\mu)$ be an increasing sequence of singular cardinals so that $\cf(\mu_i)=\theta$ for every $i\in\mu$ and $\mu=\bigcup_{i\in\mu}\mu_i$.
We may assume that $2^{\mu_i}=\mu_i^+<\mu_{i+1}$ for every $i\in\mu$.\footnote{Fix a singular cardinal $\nu>\mu$ of cofinality $\theta$. Since \textsf{GCH} holds above $\mu$, one has $2^\nu=\nu^+$. This setting reflects downward below $\mu$ to unboundedly many cardinals. Thus one can make the above mentioned assumptions on the sequence of $\mu_i$'s.}

We force with $\mathbb{Q}_{\bar{\mu}}$ followed by Magidor forcing making $\theta=\cf(\mu)$ and obtaining the assumption ${\rm tcf}(\prod_{i\in\theta}\mu_i^+,J^{\rm bd}_\theta)=\mu^+$.
If $\theta=\aleph_0$ then one can simply use Prikry forcing.
Thus, the length of the iteration should be an ordinal of cofinality $\mu^+$.
We increase $2^\mu$ to any desired point (this can be done by choosing the length of the iteration to have the desired cardinality).
Notice that $2^{\mu_i}=\mu_i^+$ remains true, as $\mathbb{Q}_{\bar{\mu}}$ is $\chi$-strategically-closed for every $\chi\in\mu$ and the component of Prikry or Magidor forcing also preserves this fact.
Thus the assumptions of Theorem \ref{thmnegfrompcf} hold in the generic extension, and the corollary follows.

\hfill \qedref{coranswertoeh5}

It seems that the above method cannot be applied to $\aleph_\omega$.
The main point is that our singular cardinal $\mu$ of cofinality $\theta$ should be a limit of singular cardinals with the same cofinality.
Thus, the negative colorings along the way are always with the same number of colors (namely, $\theta$) and hence one can produce a coloring over the cardinal $\mu^+$ with $\theta$-many colors.
Since there are no singular cardinals below $\aleph_\omega$ at all, the above proof is inapplicable as it is to this case.
However, $\aleph_{\omega^2}$ seems suitable for this pattern of proof.
Indeed, the cofinality of $\aleph_{\omega^2}$ is $\omega$ and it is a limit of singular cardinals of countable cofinality.

\begin{theorem}
  \label{thmdownto} Assuming the existence of a strong cardinal in the ground model, one can force $\mu^+\nrightarrow(\mu^+,(3)_{\cf(\mu)})^2$ with $2^\mu>\mu^+$ where $\mu$ is a strong limit cardinal, where $\mu=\aleph_{\omega^2}$.
\end{theorem}

\par\noindent\emph{Proof}. \newline
Let $\mu$ be a strong cardinal and let $\lambda\geq\mu^{++}$ be a regular cardinal.
Let $E$ be a $(\mu,\lambda)$-extender and let $\jmath:V\rightarrow{M}\cong{\rm Ult}(V,E)$ be the canonical embedding, where $M\supseteq V_{\mu^{++}}$.
We assume \textsf{GCH} in the ground model.
In order to force the above statement at $\mu$ we use the Extender-based Prikry forcing, and in order to obtain the negative relation at $\mu=\aleph_{\omega^2}$ we use the same forcing with interleaved collapses.\footnote{This is elaborated in the Appendix.}

Let $G$ be $V$-generic for this forcing notion.
Notice that $\mu$ is a strong limit cardinal in $V[G]$, and $2^\mu=\mu^{++}$.
Likewise, $\mu$ is a singular cardinal of countable cofinality in the generic extension, and \textsf{GCH} still holds below $\mu$ in $V[G]$.
We can add the collapses to make $\mu=\aleph_{\omega^2}$ in $V[G]$.

Let $(\rho_n\mid n\in\omega)$ be the Prikry sequence added through the (unique) normal ultrafilter of $E$.
It is known that ${\rm tcf}(\prod_{n\in\omega}\rho_n^{++},J^{\rm bd}_\omega)=\mu^{++}$ in the generic extension, see \cite{MR2768695}.
Moreover, up to a modification of a proper initial segment, this is the only sequence with true cofinality $\mu^{++}$ in this product.
Hence, if $(\mu_n\mid n\in\omega)$ is an increasing sequence of singular cardinals with countable cofinality such that $\mu=\bigcup_{n\in\omega}\mu_n$ then ${\rm tcf}(\prod_{n\in\omega}\mu_n^+,J^{\rm bd}_\omega)=\mu^{+}$ in $V[G]$.
For these facts we refer to \cite{yairbd}, and to the appendix at the end of the present paper.
We see that all the assumptions of Theorem \ref{thmnegfrompcf} hold, and therefore $\mu^+\nrightarrow(\mu^+,(3)_{\cf(\mu)})^2$.
In the setting of the Extender-based Prikry forcing with interleaved collapses we can make $\mu=\aleph_{\omega^2}$ in $V[G]$.
This is the first infinite cardinal which can be represented as a limit of a sequence $(\mu_n\mid n\in\omega)$ as above, so the proof is accomplished.

\hfill \qedref{thmdownto}

Notice that the large cardinal assumption here is just a strong cardinal, and in fact one needs only something around a $\mu^{++}$-strong cardinal, which is not so far from the optimal consistency strength of the failure of \textsf{SCH}.
Gitik showed in \cite{gitik} that $o(\mu)=\mu^{++}$ is sufficient for this result.
We indicate, however, that in this theorem $\cf(\mu)^{V[G]}=\omega$, and in the previous theorems $\cf(\mu)$ is arbitrary. The latter might prove more challenging, hence we started there from the assumption of supercompactness.

We conclude this section with a question about the possible consistency of the opposite situation.
Maybe the most interesting problem which issues from our study is whether the negative relation holds in \textsf{ZFC}.
We believe that the positive relation $\lambda^+\rightarrow(\lambda^+,(3)_{\cf(\lambda)})^2$ is consistent, but we do not know how to prove this:

\begin{question}
  \label{qeh5positive} Is it consistent that $\lambda$ is a strong limit singular cardinal and $\lambda^+\rightarrow(\lambda^+,(3)_{\cf(\lambda)})^2$? Is it forceable at $\lambda=\aleph_\omega$?
\end{question}

It follows from the proofs in this section that a positive answer to the above question must be forced by adding a lot of bounded subsets to $\lambda$.
In any attempt to force a positive relation one has to eliminate \textsf{GCH} at every unbounded sequence of cardinals in $\lambda$ whose true cofinality is $\lambda^+$.

\newpage

\section{More pcf and the first infinite singular cardinal}

In this section we prove that the negative relation $\lambda^+\nrightarrow(\lambda^+,(3)_{\cf(\lambda)})^2$ is consistent with $\lambda$ being a strong limit singular cardinal and $2^\lambda>\lambda^+$ even if $\lambda$ is not a limit of singular cardinals.
In particular, this setting is forceable at $\lambda=\aleph_\omega$.

The basic idea is similar to that of the proof of Theorem \ref{thmnegfrompcf} in the previous section.
Namely, from negative arrow relations on a sequence of cardinals below $\lambda$ one obtains the negative arrow relation at $\lambda^+$, provided that the true cofinality of the sequence is $\lambda^+$.
Thus, pcf theory enables us to lift combinatorial properties below $\lambda$ to $\lambda^+$.
However, in order to incorporate $\aleph_\omega$ into this framework one has to prove relevant statements about regular cardinals (e.g., the $\aleph_n$'s). For those cardinals we refine the negative partition relation to a weaker negative partition relation, relative to certain filters. Using an appropriate pcf structure that respects those filters, we will obtain the parallel of Theorem \ref{thmnegfrompcf}.

So, first we deal with obtaining those weaker instances of the negative partition relation from local instances of the generalized continuum hypothesis.
This is the content of the first claim of this section.

\begin{claim}
  \label{clmtwocolors} Assume that $\kappa<\lambda$ are regular cardinals and $\lambda=\lambda^{<\kappa}$.
  Then there exists a pair $(c,\mathscr{D})$ such that:
  \begin{enumerate}
    \item [$(a)$] $c:[\lambda]^2\rightarrow\{0,1\}$.
    \item [$(b)$] $\mathscr{D}$ is a $\kappa$-complete (proper) filter over $\lambda$.
    \item [$(c)$] If $A\subseteq\lambda$ and $c\upharpoonright[A]^2$ is constantly zero then $A=\varnothing\ \text{mod}\ \mathscr{D}$.
    \item [$(d)$] If $B=\{\alpha,\beta,\gamma\}\in[\lambda]^3$ then $0\in c''[B]^2$.
  \end{enumerate}
  Thus there is no $0$-monochromatic $\mathscr{D}$-positive set and no $1$-monochromatic triangle under $c$.
\end{claim}

\par\noindent\emph{Proof}. \newline
Enumerate the elements of $[\lambda]^{<\kappa}$ by $(u_\xi\mid\xi\in\lambda)$, where $u_0=\varnothing$ and each $u\in[\lambda]^{<\kappa}$ appears $\lambda$-many times in the enumeration.
We shall define $c$ as $\bigcup_{\alpha\in\lambda}c_\alpha$, so we define $c_\alpha$ and an ordinal $\xi_\alpha$ by induction on $\alpha\in\lambda$ as follows.
For $\alpha=0$ let $c_\alpha=\varnothing$ and if $0<\alpha$ is a limit ordinal then $c_\alpha=\bigcup\{c_\beta:\beta\in\alpha\}$.

If $\alpha=\beta+1$ then we define $c_\alpha$ and $\xi_\beta$ (so $\xi_\beta$ is picked at the $(\beta+1)$th stage).
Let $\xi_\beta$ be the minimal $\xi\in\lambda$ so that $\xi\notin\{\xi_\gamma\mid\gamma\in\beta\}\cup\{0\}$ and $c_\beta\upharpoonright[u_\xi]^2$ is constantly zero, if there is such an ordinal.
If not, let $\xi_\beta=0$.
Now if $\eta<\zeta<\beta$ then let $c_\alpha(\{\eta,\zeta\})=c_\beta(\{\eta,\zeta\})$.
For every $\gamma\in\beta$ let $c_\alpha(\{\gamma,\beta\})=1$ if $\gamma\in u_{\xi_\beta}$ and let $c_\alpha(\{\gamma,\beta\})=0$ if $\gamma\notin u_{\xi_\beta}$.
Thus $c_\alpha$ extends $c_\beta$, and the new values are determined according to the membership in $u_{\xi_\beta}$.
Let $c=\bigcup_{\alpha\in\lambda}c_\alpha$.

Let $\mathcal{I}$ be the ideal that is $\kappa$-generated by the $0$-monochromatic subsets of $\lambda$ under $c$, and the bounded subsets of $\lambda$.
Formally, for every $A\subseteq\lambda$ let $A\in\mathcal{I}$ iff there are $\zeta\in\kappa$ and $A_\varepsilon\subseteq\lambda$ for every $\varepsilon\in\zeta$ so that $c''[A_\varepsilon]^2=\{0\}$ for each $\varepsilon\in\zeta$ and $A-\bigcup\{A_\varepsilon\mid\varepsilon\in\zeta\}\in[\lambda]^{<\lambda}$.
Let $\mathscr{D}$ be the dual filter, namely $\{\lambda-A\mid A\in\mathcal{I}\}$.
Clearly, $\mathscr{D}$ is a $\kappa$-complete filter over $\lambda$.
Note that $(c)$ holds since every $0$-monochromatic set under $c$ belongs to $\mathcal{I}$ by definition.
Also, if $\varepsilon<\sigma<\tau$ and $c(\varepsilon,\tau)=c(\sigma,\tau)=1$ then $\varepsilon,\sigma\in u_{\xi_\tau}$ and then $c(\varepsilon,\sigma)=0$, thus $(d)$ holds as well.
We must prove, however, that $\mathscr{D}$ is a proper filter (or, in other words, that $\lambda\notin\mathcal{I}$).

Assume towards contradiction that $\lambda\in\mathcal{I}$.
Fix $\zeta\in\kappa$ and $A_\varepsilon\subseteq\lambda$ for every $\varepsilon\in\zeta$ such that $c\upharpoonright[A_\varepsilon]^2$ is constantly zero for each $\varepsilon\in\zeta$ and $B=\lambda-\bigcup_{\varepsilon\in\zeta}A_\varepsilon\in[\lambda]^{<\lambda}$.
Since $\lambda$ is regular, $B$ is bounded in $\lambda$.
We may assume, without loss of generality, that $|A_\varepsilon|=\lambda$ for every $\varepsilon\in\zeta$, since $B$ can be augmented by adding every $A_\varepsilon$ of size less than $\lambda$ to $B$.

We choose a sequence of ordinals $(\beta_{\alpha\varepsilon}\mid\alpha\in\lambda,\varepsilon\in\zeta)$ with no repetitions such that $\beta_{\alpha\varepsilon}\in{A_\varepsilon}$ for every $\alpha\in\lambda,\varepsilon\in\zeta$.
The choice is possible since $|A_\varepsilon|=\lambda$ for every $\varepsilon\in\zeta$.
For every $\alpha\in\lambda$ let $V_\alpha=\bigcup\{u_{\xi_{\beta_{\alpha\varepsilon}}}\mid\varepsilon\in\zeta\}$ and let $W_\alpha=\{\gamma\in\lambda\mid\exists\varepsilon\in\zeta, \beta_{\gamma\varepsilon}\in V_\alpha\}$.
Notice that $V_\alpha\in[\lambda]^{<\kappa}$ since $\kappa$ is regular, and hence $W_\alpha\in[\lambda]^{<\kappa}$ as well.
We may assume that $\alpha\notin W_\alpha$ for every $\alpha\in\lambda$ by removing this ordinal if needed.
Apply Hajnal's free subset theorem to the collection $\{W_\alpha\mid\alpha\in\lambda\}$ and let $Y\in[\lambda]^\lambda$ be free.\footnote{We use the following version: suppose that $\kappa<\lambda, f:\lambda\rightarrow[\lambda]^{<\kappa}$ is defined by $f(\alpha)=W_\alpha$ for every $\alpha\in\lambda$ and $\alpha\notin W_\alpha$ for every $\alpha\in\lambda$. Then there exists an $f$-free set $Y\in[\lambda]^\lambda$. A proof of this version appears in \cite{MR795592}.}
That is, if $\{\alpha,\beta\}\subseteq{Y}$ then $\alpha\notin{W_\beta}$.

Let $f:\zeta\rightarrow{Y}$ be an increasing function satisfying $B\cap rang(f)=\varnothing$ and, moreover, if $\varepsilon'<\varepsilon<\zeta$ then $\beta_{f(\varepsilon')\varepsilon'}<\beta_{f(\varepsilon)\varepsilon}$.
Such a function exists since $|Y|=\lambda$.
Let $u=\{\beta_{f(\varepsilon)\varepsilon}\mid\varepsilon\in\zeta\}$.
Observe that $u\in[\lambda]^{<\kappa}$ and $c\upharpoonright[u]^2$ is constantly zero.
Indeed, if $\varepsilon'<\varepsilon<\zeta$ and $c(\beta_{f(\varepsilon')\varepsilon'},\beta_{f(\varepsilon)\varepsilon})=1$ then $\beta_{f(\varepsilon')\varepsilon'}\in u_{\xi_{\beta_{f(\varepsilon)\varepsilon}}}\subseteq V_{f(\varepsilon)}$.
But then $f(\varepsilon')\in W_{f(\varepsilon)}$ by definition, and we know that $f(\varepsilon')\notin W_{f(\varepsilon)}$ as both $f(\varepsilon')$ and $f(\varepsilon)$ belong to $Y$.

Recall that $u$ appears $\lambda$ many times in our enumeration of the elements of $[\lambda]^{<\kappa}$.
At some point there will be an ordinal $\xi_\alpha$ so that $u=u_{\xi_\alpha}$ and $\xi_\alpha$ is the first ordinal for which $c\upharpoonright{[u_{\xi_\alpha}]^2}=\{0\}$, since the previous $\xi_\gamma$'s will be exhausted.
From this argument it follows that $\xi_\alpha$ can be arbitrarily large, so we choose $\xi_\alpha$ so that $\alpha>\bigcup{B}$.
Since $\alpha\notin{B}$, there is some $\varepsilon\in\zeta$ such that $\alpha\in{A_\varepsilon}$.
By definition, $\beta_{f(\varepsilon)\varepsilon}\in{A_\varepsilon}$ as well, thus $c(\{\alpha,\beta_{f(\varepsilon)\varepsilon}\})=0$ as $c\upharpoonright[A_\varepsilon]^2$ is constantly zero.
However, $\beta_{f(\varepsilon)\varepsilon}\in{u}=u_{\xi_\alpha}$, so $c(\{\alpha,\beta_{f(\varepsilon)\varepsilon}\})=1$, a contradiction.

\hfill \qedref{clmtwocolors}

The above claim is exactly what we need in order to lift negative partition relations to the negative arrow relation over a successor of a singular cardinal.
We emphasize that the coloring $c$ may possess a $0$-monochromatic set $A$ of size $\lambda$.
Only the weaker statement $A=\varnothing\ \text{mod}\ \mathscr{D}$ is guaranteed.
But this will be sufficient as shown in the following.

\begin{theorem}
  \label{thmnegrelation} Suppose that:
  \begin{enumerate}
    \item [$(A)$] \begin{enumerate}
                    \item [$(\aleph)$] $\lambda>\cf(\lambda)=\theta$.
                    \item [$(\beth)$] $\lambda=\bigcup_{i\in\theta}\kappa_i=\bigcup_{i\in\theta}\lambda_i$.
                    \item [$(\gimel)$] $\kappa_i=\cf(\kappa_i)<\lambda_i=\cf(\lambda_i)$ and $\lambda_i^{<\kappa_i}=\lambda_i$ for every $i\in\theta$.
                    \item [$(\daleth)$] $|\prod_{i<j}\lambda_i|<\kappa_j$ for every $j\in\theta$.
                  \end{enumerate}
    \item [$(B)$] \begin{enumerate}
                    \item [$(\aleph)$] $\mathscr{D}_i$ is a $\kappa_i$-complete proper filter over $\lambda_i$ for every $i\in\theta$.
                    \item [$(\beth)$] $c_i:[\lambda_i]^2\rightarrow\{0,1\}$ for every $i\in\theta$.
                    \item [$(\gimel)$] If $A\subseteq\lambda_i$ and $c_i''[A]^2=\{0\}$ then $A=\varnothing\ \text{mod}\ \mathscr{D}_i$.
                    \item [$(\daleth)$] If $\alpha<\beta<\gamma<\lambda_i$ and $\varepsilon=c_i(\{\alpha,\beta\})=c_i(\{\alpha,\gamma\}) = c_i(\{\beta,\gamma\})$ then $\varepsilon=0$.
                  \end{enumerate}
    \item [$(C)$] \begin{enumerate}
                    \item [$(\aleph)$] $\bar{\eta}=(\eta_\alpha\mid\alpha\in\lambda^+)\subseteq\prod_{i\in\theta}\lambda_i$.
                    \item [$(\beth)$] If $\alpha<\beta<\lambda^+$ then $\eta_\alpha\neq\eta_\beta$.
                    \item [$(\gimel)$] If $B_i\in\mathscr{D}_i$ for every $i\in\theta$ then there is $\alpha_0\in\lambda^+$ such that for every $\alpha_0\leq\alpha\in\lambda^+$ one can find $i_\alpha\in\theta$ so that $\eta_\alpha(i)\in{B_i}$ whenever $i_\alpha\leq{i}\in\theta$.
                  \end{enumerate}
  \end{enumerate}
  Then $\lambda^+\nrightarrow(\lambda^+,(3)_\theta)^2$.
\end{theorem}

\par\noindent\emph{Proof}. \newline
We define $c:[\lambda^+]^2\rightarrow\theta\times\{0,1\}$ as follows.
For every $\alpha<\beta<\lambda^+$ we let $i_{\alpha\beta}=\ell{g}(\eta_\alpha\cap\eta_\beta)\in\theta$, and $\varepsilon_{\alpha\beta}=c_{i_{\alpha\beta}}(\{\eta_\alpha(i),\eta_\beta(i)\})$ where $i=i_{\alpha\beta}+1$.
Now if $0<\varepsilon_{\alpha\beta}$ then we define $c(\{\alpha,\beta\})=(i_{\alpha\beta},\varepsilon_{\alpha\beta})$.
Otherwise, we set $c(\{\alpha,\beta\})=(0,0)$.
Our set of colors is $\theta\times\{0,1\}$, and $(0,0)$ stands for the first color in the negative arrow relation to be proved.

We claim that $c$ witnesses the negative arrow relation $\lambda^+\nrightarrow(\lambda^+,(3)_\theta)^2$.
To show this, we must prove two propositions.
Firstly, let us show that if $\varepsilon=1$ then the graph $(\lambda^+,\{\{\alpha,\beta\}\mid c(\{\alpha,\beta\})=(i,\varepsilon)\})$ has no triples.
Indeed, if $c(\{\alpha,\beta\})=c(\{\alpha,\gamma\})=c(\{\beta,\gamma\})=(i,\varepsilon)$ then the triple $\{\eta_\alpha(i),\eta_\beta(i),\eta_\gamma(i)\}$ is $\varepsilon$-monochromatic under $c_i$, contradicting $(B)(\daleth)$ as $\varepsilon=1$.
Secondly, we argue that if $A\in[\lambda^+]^{\lambda^+}$ then $c\upharpoonright[A]^2$ is not constantly $(0,0)$.

To see this, assume towards contradiction that $A\in[\lambda^+]^{\lambda^+}$ is a counterexample.
For every $i\in\theta$ let $A_i=\{\eta_\alpha(i)\mid\alpha\in{A}\}\subseteq\lambda_i$.
We distinguish two cases.
In the first case, $A_i\in\mathscr{D}_i^+$ for some $i\in\theta$.
Fix such $i$ and choose $\alpha_\beta\in{A}$ for every $\beta\in{A_i}$ so that $\eta_{\alpha_\beta}(i)=\beta$.
For every $\nu\in\prod_{j<i}\lambda_j$ let $A_{i\nu}=\{\beta\in{A_i}\mid\eta_{\alpha_\beta}\upharpoonright{i}=\nu\}$.
Thus $A_i=\bigcup\{A_{i\nu}\mid\nu\in\prod_{j<i}\lambda_j\}$.
Since $\mathscr{D}_i$ is $\kappa_i$-complete and $|\prod_{j<i}\lambda_j|<\kappa_i$, for some $\nu\in\prod_{j<i}\lambda_j$ one has $A_{i\nu}\in\mathscr{D}^+_i$.
By assumption $(B)(\gimel)$, there are $\beta,\beta'\in A_{i\nu}$ so that $\beta\neq\beta'$ and $c_i(\{\beta,\beta'\})=\varepsilon>0$.
By definition, $c(\{\alpha_\beta,\alpha_{\beta'}\})=(i,\varepsilon)\neq(0,0)$.
But $\alpha_\beta,\alpha_{\beta'}\in{A}$, so $c(\{\alpha_\beta,\alpha_{\beta'}\})=(0,0)$, a contradiction.

In the second case, $A_i=\varnothing\ \text{mod}\ \mathscr{D}_i$ for every $i\in\theta$.
Recall that $\forall^*i<\theta$ means for every $i\in\theta$ apart from boundedly many ordinals of $\theta$.
Define $B=\{\alpha\in\lambda\mid(\forall^{*}i<\theta)(\eta_\alpha(i)\in{A_i})\}$.
From $(C)(\gimel)$ we see that $B\in[\lambda]^{<\lambda}$, as $B_i=\lambda_i-A_i\in\mathscr{D}_i$ for every $i\in\theta$ according to the second case.
Now if $\alpha\in{A}$ then, by definition, $\eta_\alpha(i)\in{A_i}$ for every $i\in\theta$.
Hence $A\subseteq{B}$ and consequently $|A|<\lambda$, a contradiction.

\hfill \qedref{thmnegrelation}

In order to utilize the above theorem one has to show that the assumptions there hold (or can be forced) under the relevant circumstances.
The assumptions of $(A)$ are easily satisfied for an appropriate choice of cardinals when $\lambda$ is a strong limit singular cardinal.
The assumptions of $(B)$ follow from Claim \ref{clmtwocolors}.
Thus the only challenge is $(C)$, and this is our next goal.

\begin{lemma}
  \label{lembasicomponent} Let $\lambda$ be a strongly inaccessible cardinal.
  Let $R\subseteq\lambda$ be an unbounded set of regular cardinals.
  For every $\zeta\in{R}$ let $\mathscr{D}_\zeta$ be a $\zeta$-complete uniform filter over $\zeta^+$.
  There is a $(<\lambda)$-strategically-closed $\lambda^+$-cc forcing notion $\mathbb{D}$, for which the following hold in the generic extension by $\mathbb{D}$:
  \begin{enumerate}
    \item [$(\aleph)$] $(f_\alpha\mid\alpha\in\lambda^+)$ is a scale in $\prod_{\zeta\in{R}}\zeta^+$.
    \item [$(\beth)$] For every $(A_\zeta\mid\zeta\in{R})\in\prod_{\zeta\in{R}}\mathscr{D}_\zeta$ there is $\alpha_0\in\lambda^+$ such that for each $\alpha_0\leq\alpha\in\lambda^+$ there exists $i_\alpha\in\lambda$ such that if $i_\alpha\leq\zeta\in{R}$ then $f_\alpha(\zeta)\in{A_\zeta}$.
  \end{enumerate}
\end{lemma}

\par\noindent\emph{Proof}. \newline
Let $(\mathbb{D}_\alpha,\name{\mathbb{H}}_\beta\mid\beta\in\lambda^+,\alpha\leq\lambda^+)$ be a $(<\lambda)$-support iteration of length $\lambda^+$, where $\name{\mathbb{H}}_\beta$ is (a name of) a forcing notion defined in the generic extension by $\mathbb{D}_\beta$ as follows.
A condition $p\in\mathbb{H}_\beta$ is a pair $(s,\bar{A})=(s_p,\bar{A}_p)$, where $s$ is a function with ${\rm dom}(s)=R\cap\xi_s$ for some $\xi_s\in\lambda$ and $\bar{A}$ is a sequence of sets of the form $(A_\zeta\mid\zeta\in{R}-{\rm dom}(s))$.
For every $\zeta\in{\rm dom}(s)$ one has $s(\zeta)\in\zeta^+$ and for every $\zeta\in{R}-{\rm dom}(s)$ one has $A_\zeta\in\mathscr{D}_\zeta$.
If $p,q\in\mathbb{H}_\beta$ then $p\leq{q}$ iff $s_p\trianglelefteq{s_q}, A^q_\zeta\subseteq A^p_\zeta$ whenever $\zeta\in{R}-{\rm dom}(s_q)$ and if $\zeta\in{\rm dom}(s_q)-{\rm dom}(s_p)$ then $s_q(\zeta)\in A^p_\zeta$.
Thus $\mathbb{H}_\beta$ approximates a function in $\prod_{\zeta\in{R}}\zeta^+$ in a Hechlerish style.
We shall say that $s_p$ is \emph{the stem} of $p$ and the sequence $\bar{A}_p$ is \emph{the pure part} of $p$.

Observe that $\name{\mathbb{H}}_\beta$ is (forced to be) $(<\lambda)$-strategically-closed.
In fact, if $\gamma\in\lambda$ then there exists a $\gamma^+$-directed-closed dense open set of conditions, that is, $T_\gamma=\{p\in\mathbb{H}_\beta\mid\bigcup{\rm dom}(s_p)>\gamma\}$.
To see this, suppose that $\{p_i\mid i<\gamma\}$ is directed.
Define $s_p=\bigcup_{i<\gamma}s_{p_i}$ and for every $\zeta\in{R}-{\rm dom}(s_p)$ let $A^p_\zeta=\bigcap_{i<\gamma}A^{p_i}_\zeta$.
Notice that $A^p_\zeta\in\mathscr{D}_\zeta$ since $\zeta>\gamma$ and $\mathscr{D}_\zeta$ is $\zeta$-complete.
Now $p=(s_p,\bar{A}_p)$ is an upper bound for $\{p_i\mid i<\gamma\}$.
It is easy to see that $\mathbb{H}_\beta$ is $\lambda$-centered.
Indeed, if $p,q\in\mathbb{H}_\beta$ and $s_p=s_q$ then $p\parallel{q}$.
Since the number of stems is $\lambda$ (recall that $\lambda$ is strongly inaccessible) we conclude that $\mathbb{H}_\beta$ is $\lambda$-centered.

Let $H\subseteq\mathbb{H}_\beta$ be generic.
Define $f_\beta=\bigcup\{s\mid(s,\bar{A})\in{H}\}$.
By the directness of $H$ and simple density arguments, $f_\beta$ is a function, ${\rm dom}(f_\beta)=R$ and $f_\beta(\zeta)\in\zeta^+$ for every $\zeta\in{R}$.
Let $\mathbb{D}=\mathbb{D}_{\lambda^+}$.
Observe that $\mathbb{D}$ is $(<\lambda)$-strategically-closed and $\lambda^+$-cc.
Indeed, $\mathbb{D}$ is a $(<\lambda)$-support iteration and each component in the iteration is $(<\lambda)$-strategically-closed and $\lambda$-centered.
Fix a $V$-generic set $G\subseteq\mathbb{D}$.
Let us show that the statement of the lemma holds in $V[G]$.

As noted above, $(f_\alpha\mid\alpha\in\lambda^+)\subseteq\prod_{\zeta\in{R}}\zeta^+$.
Since $\mathbb{D}$ is $(<\lambda)$-strategically-closed, no new bounded subsets of $\lambda$ are introduced in $V[G]$.
However, new \emph{sequences} of measure-one sets (of length $\lambda$) are introduced.
Fix a sequence $\name{\bar{A}}=(\name{\bar{A}}_\zeta\mid\zeta\in{R})$ such that $\name{\bar{A}}[G]\in\prod_{\zeta\in{R}}\mathscr{D}_\zeta\cap V[G]$.
Since $\mathbb{D}$ is $\lambda^+$-cc, there must be some $\gamma\in\lambda^+$ so that $\name{\bar{A}}\in V^{\mathbb{D}_\gamma}$.
We claim that if $\gamma\leq\beta\in\lambda^+$ then there exists $i_\beta\in\lambda$ such that $f_\beta(\zeta)\in\name{A}_\zeta$ whenever $i_\beta\leq\zeta\in\lambda$.
If we prove this statement then the proof of the lemma will be accomplished.

To see this, let $p$ be an arbitrary condition that forces $\name{\bar{A}}$ in the $\gamma$th stage of the iteration to be a sequence of measure-one sets.
Fix $\beta\in[\gamma,\lambda^+)$.
Define $q\in\mathbb{D}$ so that $p\leq{q}, \beta\in{\rm dom}(q)$ and let $q(\beta)=(t_q,\name{\bar{B}}_q)$.
We require that $q\upharpoonright\beta\Vdash \name{B}^q_\xi\subseteq\name{A}_\xi$ for every $\xi\in{\rm dom}(\name{\bar{B}}_q)$.
This is possible as $\name{A}_\xi$ is a $\mathbb{D}_\gamma$-name and hence also a $\mathbb{D}_\beta$-name.
Finally, let $i_\beta=\bigcup{\rm dom}(t_q)$.
If $i_\beta\leq\zeta\in\lambda$ then $q\Vdash f_\beta(\zeta)\in \name{B}^q_\zeta\subseteq\name{A}_\zeta$, so as $p$ was arbitrary we are done.

\hfill \qedref{lembasicomponent}

We remark that for densely many conditions $p\in\mathbb{D}$, for every $\alpha\in{\rm dom}(p)$, the stem of the condition $p(\alpha)$ is a canonical name of a ground model function.
Hence we may assume, without loss of generality, that this is true for every condition in $\mathbb{D}$.

By the previous lemma, we can force our pcf assumption using a $(<\lambda)$-strategically-closed forcing notion.
We would like, therefore, to prove that a supercompact cardinal $\lambda$ will remain supercompact after such a forcing.
This will be instrumental since we will singularize this cardinal with various kinds of Prikry-type forcings.
Thus we have to show that a supercompact cardinal $\lambda$ can be forced to be indestructible under such forcings.
We shall follow the ideas of Laver, with slight changes, in order to force this property over a supercompact cardinal in the ground model.

As our forcing is not sufficiently directed-closed, we need a weaker property that would imply the existence of a master condition for lifting certain supercompact embeddings.
\begin{lemma}
  \label{lemaster} Let $R$ and $\mathbb{D}$ be as above, and let $\zeta\in{R}$.
  Suppose that $F\subseteq\mathbb{D}$ is a directed set of conditions, $|F|<\zeta$, and for every $\alpha\in\bigcup\{{\rm dom}(p)\mid p\in{F}\}$ one has $\bigcup\{{\rm dom}(s_{p(\alpha)})\mid p\in{F}, \alpha\in {\rm dom}(p)\}\geq\zeta$.
  Then there exists a master condition $m\in\mathbb{D}$ so that $p\leq{m}$ for every $p\in{F}$.
\end{lemma}

\par\noindent\emph{Proof}. \newline
We commence with letting ${\rm dom}(m)=\bigcup\{{\rm dom}(p)\mid p\in{F}\}$.
For every $\alpha\in{\rm dom}(m)$ let $s_{m(\alpha)}=\bigcup\{s_{p(\alpha)}\mid p\in{F}, \alpha\in {\rm dom}(p)\}$ and let $\bar{A}^{m(\alpha)}_\xi=\bigcap\{\bar{A}^{p(\alpha)}_\xi\mid p\in{F}, \alpha\in{\rm dom}(p)\}$, for every $\xi\in R-\bigcup {\rm dom}(s_{m(\alpha)})$.
Notice that $\bar{A}^{m(\alpha)}_\xi$ is a measure-one set in $\mathscr{D}_\xi$, by the completeness of $\mathscr{D}_\xi$.
We prove by induction on $\delta\in\lambda^+$ that $m\upharpoonright\delta$ is a condition and $p\upharpoonright\delta\leq m\upharpoonright\delta$ for every $p\in{F}$.

For $\delta=0$ this is trivial, and for a limit ordinal $\delta$ it follows from the fact that $|{\rm dom}(m)|<\lambda$.
We are left with the successor steps, so fix $\delta=\alpha+1\in\lambda^+$.
By the induction hypothesis, $m\upharpoonright\alpha$ forces $\bar{A}_{p(\alpha)}$ to be a sequence of measure-one sets for every $p\in{F}$.
Hence $m\upharpoonright\delta$ is a condition.
It remains to show that $p\upharpoonright\delta\leq m\upharpoonright\delta$ for every $p\in{F}$.

Fix $p\in{F}$.
From the definition it follows that $s_{p(\alpha)}\trianglelefteq s_{m(\alpha)}$ and $\bar{A}^{m(\alpha)}_\xi\subseteq \bar{A}^{p(\alpha)}_\xi$ for every relevant $\xi$.
Suppose that $\xi\in{\rm dom}(s_{m(\alpha)})-{\rm dom}(s_{p(\alpha)})$.
Since $s_{m(\alpha)}$ is the union of $s_{p(\alpha)}$s, there exists $r\in{F}$ for which $\xi\in{\rm dom}(s_{r(\alpha)})$ and $s_{m(\alpha)}(\xi)=s_{r(\alpha)}(\xi)$.
Since $F$ is directed, there exists $t\in{F}$ such that $r,p\leq{t}$.
It follows that $s_{t(\alpha)}(\xi)=s_{r(\alpha)}(\xi)$, and hence $t\upharpoonright\alpha\Vdash{m(\alpha)}(\xi)=s_{t(\alpha)}(\xi)\in \bar{A}^{p(\alpha)}_\xi$, so we are done.

\hfill \qedref{lemaster}

As mentioned before, we need the above lemma in order to prove a version of Laver's indestructibility.
In the work of Laver, a supercompact cardinal $\lambda$ is forced to be indestructible under any $\lambda$-directed-closed forcing notion.
One cannot expect $\lambda$ to be indestructible under every $\lambda$-strategically-closed forcing, since one can force a non-reflecting stationary set with such a forcing notion.
However, if the set $R$ is sufficiently sparse then one can prepare the universe so that the corresponding forcing notion $\mathbb{D}$ will preserve the supercompactness of $\lambda$.
This is the content of our next lemma, and the main idea is that we define the Laver iteration only with respect to some strategically closed forcing notions.

\begin{lemma}
  \label{lemlaver} Let $\lambda$ be supercompact and assume $\mathsf{GCH}$.
  There is a forcing notion $\mathbb{L}$ such that the following holds in the generic extension by $\mathbb{L}$.
  If $R\subseteq\lambda$ is a set of double-double successors\footnote{For the lemma as it is, double successors are sufficient; but we shall apply the lemma to a set of double-double successors, so we phrase it accordingly.} of strongly inaccessible cardinals and $\mathscr{D}_\zeta$ is a $\zeta$-complete filter over $\zeta^+$ for every $\zeta\in{R}$ then the associated forcing $\mathbb{D}$ preserves the supercompactness of $\lambda$ as well as $\mathsf{GCH}$.
\end{lemma}

\par\noindent\emph{Proof}. \newline
Let $\ell:\lambda\rightarrow{V_\lambda}$ be a Laver function.\footnote{See \cite[Lemma, p. 386]{MR0472529}.}
Let $\mathbb{L}$ be the Easton support iteration $(\mathbb{L}_\alpha,\name{\mathbb{Q}}_\beta\mid\beta<\lambda,\alpha\leq\lambda)$ where $\name{\mathbb{Q}}_\beta$ is an $\mathbb{L}_\beta$-name of the trivial forcing unless $\beta$ is strongly inaccessible, $\ell(\alpha)\in V_\beta$ for each $\alpha\in\beta$ and $\ell(\beta)$ is a pair of the form $(\gamma,\name{\tau})$ and $\name{\tau}$ is an $\mathbb{L}_\beta$-name of some $(<\beta)$-strategically-closed forcing notion that preserves \textsf{GCH}.

One can verify that $\mathbb{L}$ is $\lambda$-cc and if $\beta$ is an inaccessible closure point of $\ell$ then $\beta$ is preserved by $\mathbb{L}$.
It follows from the choice of the $\name{\mathbb{Q}}_\beta$s and standard arguments that \textsf{GCH} holds in the generic extension by $\mathbb{L}$.
Choose a $V$-generic set $G\subseteq\mathbb{L}$, and in $V[G]$ choose a $V[G]$-generic set $g\subseteq\mathbb{D}$.
We claim that $\lambda$ is supercompact in $V[G][g]$.

To see this, suppose that $\mu$ is an arbitrary ordinal and $\jmath:V\rightarrow{M}$ is a $\mu$-supercompact elementary embedding, with $\jmath(\ell)(\lambda)=(\mu,\mathbb{D})$.
By increasing $\mu$ if needed we may assume that $\mu=\cf(\mu)\geq\lambda^+$.
We may also assume that $\jmath$ is an ultrapower embedding derived from some fine and normal measure over $\mathcal{P}_\lambda\mu$.

By elementarity, $\jmath(\mathbb{L})=\mathbb{L}\ast\mathbb{D}\ast\mathbb{L}_{tail}$, where $\mathbb{L}_{tail}$ is $\mu^+$-strategically-closed in $V[G][g]$.
Let $F=\{\jmath(\name{q})\mid\name{q}\in{g}\}$.
By the closure of $M$ and the fact that $\jmath\upharpoonright\name{\mathbb{D}}\in{M}$ we see that $F\in M[G][g]$.
Moreover, $F$ satisfies the assumptions of Lemma \ref{lemaster}.
Indeed, by the choice of $R$ we know that $\min(R-\lambda)>\lambda^+$ and $|F|=\lambda^+$.
Hence Lemma \ref{lemaster} applies and there is a master condition $m$ for $F$.

By \textsf{GCH}, which holds in $V[G][g]$, the number of dense subsets of $\mathbb{L}_{tail}$ is $|\jmath(2^{\lambda^+})|^V=|{}^{\mathcal{P}_\lambda\mu}(\lambda^{++})|^V=\mu^+$.
Using the fact that $\mathbb{L}_{tail}$ is $\mu^+$-closed, one can construct an $M[G][g]$-generic set $h$ containing $m$.
By the nature of $m$, $\jmath''g$ is contained in $h$.
Hence Silver's criterion applies and $\jmath:V\rightarrow{M}$ lifts to $\jmath^+:V[G][g]\rightarrow M[H][h]$.
Thus $\lambda$ is $\mu$-supercompact in $V[G][g]$.
But $\mu$ was arbitrary, hence one concludes that $\lambda$ is supercompact in $V[G][g]$ and the proof is accomplished.

\hfill \qedref{lemlaver}

The above setting is rendered in a universe satisfying \textsf{GCH}.
Our goal, however, is to force the negative arrow relation when $2^\lambda>\lambda^+$.
Furthermore, we wish to force this relation at $\lambda=\aleph_\omega$.
We shall use the extender-based Prikry forcing in order to increase $2^\lambda$, and the corresponding version with interleaved collapses in order to push down the result to $\aleph_\omega$.
The presentation of these forcing notions is inspired by the work of Merimovich, see e.g.\  \cite{MR4264147}.

Let $W$ be the universe obtained by the preparatory forcing $\mathbb{L}$.
Fix $R\subseteq\lambda$ such that every element of $R$ is of the form $\rho^{+4}$, where $\rho$ is a strongly inaccessible cardinal, and $R$ is unbounded in $\lambda$.
For every $\rho^{+4}\in{R}$ let $c_\rho:[\rho^{+4}]^2\rightarrow\{0,1\}$ be as guaranteed in Claim \ref{clmtwocolors}.
Namely, $c_\rho$ has no $1$-monochromatic triangle and if $A\subseteq\rho^{+4}$ is $0$-monochromatic under $c_\rho$ then $A=\varnothing\ \text{mod}\ \mathscr{D}_\rho$, where $\mathscr{D}_\rho$ is a $\rho^{+3}$-complete proper filter over $\rho^{+4}$.
We denote the dual (proper) ideal by $\mathcal{I}_\rho$, so $\mathcal{I}_\rho$ is $\rho^{+3}$-generated from the bounded subsets of $\rho^{+4}$ and the $0$-monochromatic sets of the coloring $c_\rho$.

Let $\mathbb{D}$ be the forcing notion associated with $R$, and let $V$ be the extension of $W$ by $\mathbb{D}$.
Let $(f_\alpha\mid\alpha\in\lambda^+)$ be the generic scale.
Fix a $(\lambda,\lambda^{++})$-extender $E$ in $V$ and let $\jmath:V\rightarrow{M}$ be the extender ultrapower map.
Let $\mathscr{U}_0$ be the (unique) normal measure of $E$, and let $\imath:V\rightarrow Ult(V,\mathscr{U}_0)\cong{N}$ be the corresponding ultrapower embedding.
Let $K_0$ be $N$-generic for the forcing $Col^N((\lambda^{+5})^N,<\imath(\lambda))$.
Let $k:N\rightarrow{M}$ be the quotient map, and let $K$ be the generic set generated by $k[K_0]$.

$$
\xymatrix{
{\rm \bf V} \ar[ddr]_{\imath} \ar[r]^\jmath & M \\ \\
& {Ult(V,\mathscr{U}_0)\cong N} \ar[uu]_{k} }
$$

\begin{claim}
  \label{clmmgeneric} $K$ is $M$-generic for the forcing $Col^M((\lambda^{+5})^M,<\jmath(\lambda))$.
\end{claim}

\par\noindent\emph{Proof}. \newline
Let $D$ be a dense subset of $Col^M((\lambda^{+5})^M,<\jmath(\lambda))$, and assume that $D\in{M}$.
We must show that $D\cap{K}\neq\varnothing$.
By the nature of $\jmath$, there are $a\in[\lambda^{++}]^n$ for some $n\in\omega$ and $f:[\lambda]^n\rightarrow{V}$ so that $\jmath{f}(a)=D$.
Define $D'=\bigcap\{\imath{f}(z)\mid z\in[(\lambda^{++})^N]^n, \imath{f}(z)\ \text{is a dense open subset of}\ Col^N((\lambda^{+5})^N,<\imath(\lambda))\}$.

Clearly $D'\in{N}$, and since $Col^N((\lambda^{+5})^N,<\imath(\lambda))$ is $(\lambda^{+5})^N$-distributive and $D'$ is the intersection of less than $(\lambda^{+})^N$ dense open sets, we see that $D'$ is dense and open.
Thus $D'\cap{K}\neq\varnothing$ and hence $k(D')\cap{K}\neq\varnothing$.
Notice that $D=k(\imath{f})(a)$ for $a\in k([(\lambda^{++})^N]^n)$, hence $D\supseteq k(D')$.
In particular, $D\cap{K}\neq\varnothing$ as sought.

\hfill \qedref{clmmgeneric}

Working in $V$, let $\mathbb{P}$ be the extender-based Prikry forcing with interleaved collapses, using the extender $E$ and the set $K$ as a guiding generic.
Following \cite{MR4264147}, for every $d\in[\lambda^{++}]^{\leq\lambda}$ we let $E(d)=\{x\in V_\lambda\mid (\jmath\upharpoonright{d})^{-1}\in\jmath(x)\}$.
Thus $E(d)$ is a $\lambda$-complete ultrafilter over $V_\lambda$, and $(\jmath\upharpoonright{d})^{-1}$ is a partial increasing function from $\jmath(d)$ to $d$.
Since $E(d)$ concentrates on such objects we may assume, without loss of generality, that if $A\in E(d)$ then $A$ is a set of partial increasing functions with domain contained in $d$.

We define our forcing notion $\mathbb{P}$ as follows.
A condition $p\in\mathbb{P}$ is a quadruple of the form $(\bar{c}^p,f^p,A^p,F^p)$, where:
\begin{enumerate}
  \item [$(a)$] $f^p$ is a partial function from $\lambda^{++}-\lambda$ to $\lambda^{<\omega}$.
  \item [$(b)$] $d^p={\rm dom}(f^p)$ satisfies $|d^p|\leq\lambda$.
  \item [$(c)$] $\lambda\in d^p$ for every $p$ and we let $f^p(\lambda)=(\rho^p_i\mid i\in{n})$.
  \item [$(d)$] $\bar{c}^p\in\prod_{i\leq{n}}Col((\rho^p_{i-1})^{+5},<\rho^p_i)$, where $\rho^p_{-1}=\omega$ and $\rho^p_n=\lambda$.
  \item [$(e)$] $A^p\in E(d^p)$.
  \item [$(f)$] If $\eta\in A^p$ then $\eta(\lambda)>\bigcup{\rm dom}(c^p_n)$, and $\eta(\lambda)>\rho^p_{n-1}$.
  \item [$(g)$] $F^p$ is a function from $A^p$ to $V_\lambda$.
  \item [$(h)$] For every $\eta\in A^p$, $F^p(\eta)\in Col(\eta(\lambda)^{+5},<\lambda)$.
  \item [$(i)$] $\jmath(F^p)((\jmath\upharpoonright{d^p})^{-1})\in{K}$.
\end{enumerate}
We have to define the forcing order.
To this end, we define the direct order $\leq^*$ and one point extensions.
Then we let the forcing order be a finite sequence of these two extensions.

Assume, therefore, that $p,q\in\mathbb{P}$.
We shall say that $q$ is a direct extension of $p$ (denoted by $p\leq^*q$) iff the following hold:
\begin{enumerate}
  \item [$(\alpha)$] $f^p\subseteq f^q$. In particular, $\ell{g}(f^p(\lambda))=\ell{g}(f^q(\lambda))$.
  \item [$(\beta)$] $A^q\subseteq\{\eta\mid\eta\upharpoonright d^p\in A^p\}$.
  \item [$(\gamma)$] $c^q_i$ extends $c^p_i$ for every ${i}\leq{n}$.\footnote{We include the extension of the collapses in the definition of the direct order.}
  \item [$(\delta)$] $F^q(\eta)\geq F^p(\eta\upharpoonright d^p)$ for every $\eta\in A^q$ (here $\geq$ is the order of the collapse).
\end{enumerate}
Assume now that $p\in\mathbb{P}$ and $\eta\in A^p$.
We define the one-point extension (or Prikry extension) $p^\frown\eta$ as a condition $r$ satisfying the following requirements:
\begin{enumerate}
  \item [$(\alpha)$] $f^r(\alpha)=f^p(\alpha)$ if $\alpha\notin{\rm dom}(\eta)$, and $f^r(\alpha)=f^p(\alpha)^\frown\eta(\alpha)$ if $\alpha\in{\rm dom}(\eta)$.
  \item [$(\beta)$] $\bar{c}^r={\bar{c}^{p}} {}^\frown F^p(\eta)$.
  \item [$(\gamma)$] $A^r=A^p-\{\rho\in A^p\mid\rho(\lambda)\leq\max\{\eta(\lambda),\bigcup{\rm dom}(F^p(\eta))\}\}$.
  \item [$(\delta)$] $F^r=F^p\upharpoonright A^r$.
\end{enumerate}
Finally, $q$ extends $p$ iff $q$ can be obtained from $p$ by a finite sequence of direct extensions and one-point extensions.
Clearly, if $p\leq_{\mathbb{P}}q$ then the process of deriving $q$ from $p$ can be presented as one single direct extension followed by a finite sequence of one-point extensions.

Notice that $(\mathbb{P},\leq^*)$ is $\sigma$-closed since $E(d)$ is $\lambda$-complete and the collapses are $\aleph_1$-complete.
It is known that $\mathbb{P}$ satisfies the strong Prikry property, preserves cardinals above $\lambda$ (including $\lambda$ itself), forces $\lambda=\aleph_\omega$, preserves \textsf{GCH} below $\lambda$ and forces $2^{\aleph_\omega}=\aleph_{\omega+2}$.

\begin{claim}
  \label{clmspp} Let $\name{\bar{\rho}}=(\rho_n\mid n\in\omega)$ be a name for the normal Prikry sequence.
  Suppose that $p\in\mathbb{P}$ and $\name{\tau}$ is a $\mathbb{P}$-name of an ordinal in $\rho^{+5}_n$.
  There exists a condition $q\in\mathbb{P}$ so that $p\leq^*q$ and a function $g$ with ${\rm dom}(g)=\lambda$, such that $g(\alpha)\in[\alpha^{+5}]^{\leq\alpha^{++}}$ for every $\alpha\in\lambda$ and $q\Vdash\name{\tau}\in g(\rho_n)$.
\end{claim}

\par\noindent\emph{Proof}. \newline
Let $D$ be a dense open set deciding the value of $\name{\tau}$.
By the strong Prikry property there are a natural number $n$ and a condition $q$ such that $p\leq^*q$ and every extension of $q$ by $n$-many one-point extensions belongs to $D$.
Recall that $d^p$ is the domain of $f^p$.
Let $\eta_0,\ldots,\eta_{n-1}$ be a sequence of $n$-many one-point extensions of $q$ such that the value of $\name{\tau}$ is decided by $q$ and these $\eta_i$'s.

Since $d^p\in[\lambda^{++}]^{\leq\lambda}$ and $((\jmath\upharpoonright{d^p})^{-1})\in\jmath(A^q)$ we may assume (by shrinking $A^q$ if needed) that for every $\eta\in{A^q}$ and every $\alpha\in{\rm dom}(\eta)\subseteq{d^p}$ it is true that $\eta(\alpha)\leq\eta(\lambda)^{++}$.
Moreover, by fixing a bijection $h:\lambda\rightarrow{d^p}$ and shrinking $A^q$ further we may assume that ${\rm dom}(\eta)=h''\eta(\lambda)$.
Consequently, $|\{\eta\in{A^q}\mid\alpha=\eta(\lambda)\}|\leq|[\alpha^{++}]^\alpha|=\alpha^{++}$ for every $\alpha\in\lambda$.
Hence we can define $g(\alpha)$ as the set of all possible values for $\name{\tau}$ using $\eta_0,\ldots,\eta_{n-1}$ with $\eta_{n-1}(\lambda)=\alpha$, and get $q\Vdash\name{\tau}\in g(\rho_n)$ as required.

\hfill \qedref{clmspp}

For every $n\in\omega$ let $\mathscr{D}_n$ be the $\rho^{+3}_n$-complete filter over $\rho^{+4}_n$ derived from Claim \ref{clmtwocolors}, and let $\mathcal{I}_n$ be the dual ideal of $\mathscr{D}_n$.
The forcing $\mathbb{P}$ certainly introduces new bounded subsets of $\lambda$, so we must show that the relevant assumptions of Theorem \ref{thmnegrelation} remain true in the generic extension by $\mathbb{P}$.
We commence with the following:

\begin{claim}
  \label{clmecoloring} Let $c:[\rho^{+4}_n]^2\rightarrow\{0,1\}$ be the coloring described in Claim \ref{clmtwocolors}.
  Let $\name{I}$ be a name of a $0$-monochromatic set under $c$.
  Then $\name{I}$ is (forced to be) contained in the union of less than $\rho^{+3}_n$-many $0$-monochromatic sets from the ground model.
  Moreover, $\name{I}\subseteq{I}$ for some $I\in{V}$ so that $I=\varnothing\ \text{mod}\ \mathscr{D}$, where $\mathscr{D}$ is the filter constructed in Claim \ref{clmtwocolors}.
\end{claim}

\par\noindent\emph{Proof}. \newline
Fix an enumeration of $\mathcal{P}^V(\rho^{+4}_n)$, for every $n\in\omega$, in the ground model.
We are assuming \textsf{GCH} in the ground model, hence the size of the enumeration is $\rho^{+5}_n$.
Let $p$ be an arbitrary condition that forces $\name{I}$ to be $0$-monochromatic.
By extending $f^p(\lambda)$ if needed we may assume that $\ell{g}(f^p(\lambda))>n$.
In other words, $\rho_n$ is decided by $p$.

Let $\mathbb{C}=\prod_{i\leq{n}}Col(\rho^{+4}_{i-1},<\rho_i)$, that is, $\mathbb{C}$ is the product of the first $n+1$ collapses mentioned in the condition $p$. Note that $|\mathbb{C}| = \rho_n$ and the next collapses are at least $\rho_n^{+5}$-closed.
Let $((s_\alpha,\xi_\alpha)\mid\alpha\in\rho^{+4}_n)$ enumerate $\mathbb{C}\times\rho^{+4}_n$.
By induction on $\alpha\in\rho^{+4}_n$ we define a condition $p_\alpha$ such that:
\begin{enumerate}
  \item [$(\aleph)$] $(p_\alpha\mid\alpha\in\rho^{+4}_n)$ is $\leq^*$-increasing.
  \item [$(\beth)$] $\bar{c}^{p_\alpha}\upharpoonright{n+1}$ is constant.
  \item [$(\gimel)$] If there is a direct extension $r$ of $p_\alpha$ which is the same as $p_\alpha$ except that $\bar{c}^{r}\upharpoonright{n+2}=s_\alpha^\frown c^{p_\alpha}_{n+1}$ and $r$ decides the truth value of the statement $\check{\xi}_\alpha\in\name{I}$ then $p_\alpha$ already decides this statement.
\end{enumerate}
Let $q$ be an upper bound of $(p_\alpha\mid\alpha\in\rho^{+4}_n)$.
Such an upper bound exists since all the relevant components are sufficiently closed.
Let $q[s]$ denote the condition obtained from $q$ upon replacing $\bar{c}^q$ by $s^\frown(\bar{c}^q\upharpoonright[n,\ell{g}(\bar{c}^q)))$.

Let $r_0 = \bar{c}^q\upharpoonright{n+1}$.
For every $s\in\mathbb{C}$ let $I_s=\{\alpha\in\rho^{+4}_n\mid q[s]\Vdash\check{\alpha}\in\name{I}\}$, so $I_s\in{V}$ and $\Vdash_{\mathbb{P}}\name{I}\subseteq\bigcup\{I_s\mid r_0\leq_{\mathbb{C}}s\}$.
Furthermore, $I_s$ is $0$-monochromatic under $c$ for every $s$ such that $r_0\leq_{\mathbb{C}}s$, hence $I_s\in\mathcal{I}_n$, where $\mathcal{I}_n$ is the dual of $\mathscr{D}_n$.
We conclude, therefore, that $\name{I}$ is (forced to be) covered by a ground model set in $\mathcal{I}_n$, as required.

\hfill \qedref{clmecoloring}

Let $R=\{\rho^{+4}_n\mid n\in\omega\}$.
We claim that $(f_\alpha\upharpoonright{R}\mid\alpha\in\lambda^+)$ is a scale satisfying the required assumption of Theorem \ref{thmnegrelation}.
Specifically, we must show that for every sequence of sets $(A_n\mid n\in\omega)$ where $A_n\in\mathscr{D}_n$ for each $n\in\omega$, there exists $\gamma\in\lambda^+$ such that if $\gamma\leq\beta\in\lambda^+$ then for some $n_\beta\in\omega$ one has $f_\beta(n)\in{A_n}$ whenever $n_\beta\leq{n}\in\omega$.

By Claim \ref{clmecoloring} we may assume that $A_n\in{V}$ for every $n\in\omega$.
Notice, however, that the whole sequence $(A_n\mid n\in\omega)$ need not be in $V$.
Let $p$ be an arbitrary condition in $\mathbb{P}$.
Applying Claim \ref{clmspp} repeatedly and using the enumeration of the ground model sets, we construct a $\leq^*$-increasing sequence $(q_n\mid n\in\omega)$ such that $q_0\geq{p}$ and $q_n\Vdash \name{A}_n\in\mathcal{A}_n(\name{\rho}_n)$, where $|\mathcal{A}_n(\gamma)|\leq\gamma^{++}$ for every $\gamma\in\lambda$.

Since $(\mathbb{P},\leq^*)$ is $\sigma$-closed, there exists a single condition $q$ such that $q_n\leq{q}$ for each $n\in\omega$.
For every $\zeta\in{R}$ let $B_\zeta=\bigcap\{E_n\mid n\in\omega, E_n\in\mathcal{A}_n(\zeta)\}$.
By the closure of $\mathscr{D}_\zeta$ we see that $B_\zeta\in\mathscr{D}_\zeta$.
By the construction, $(B_\zeta\mid\zeta\in{R})$ belongs to the ground model and hence there is an ordinal $\gamma\in\lambda^+$ such that for every $\gamma\leq\beta\in\lambda^+$ there exists $\zeta_\beta\in\lambda$ so that $f_\beta(\zeta)\in B_\zeta$ whenever $\zeta_\beta\leq\zeta\in{R}$.
This completes the proof, as $B_{\rho^{+4}_n}\subseteq{A_n}$ for every $n\in\omega$.

In order to force the failure of \textsf{SCH} at $\lambda$, as done in our results, one has to assume the existence of a measurable cardinal $\kappa$ with $o(\kappa)=\kappa^{++}$ in the ground model.
This fundamental result was proved by Gitik in \cite{MR1007865} and in \cite{MR1098782}.
In our constructions we started from a supercompact cardinal in the ground model.
The gap between these large cardinals invites the following:

\begin{question}
  \label{qstrength} Let $\lambda$ be a strong limit singular cardinal.
  \begin{enumerate}
    \item [$(\aleph)$] What is the consistency strength of the negative arrow relation $\lambda^+\nrightarrow(\lambda^+,(3)_{\cf(\lambda)})^2$ with $2^\lambda>\lambda^+$?
    \item [$(\beth)$] What is the consistency strength of the same negative relation with $2^\lambda>\lambda^+$ where $\lambda=\aleph_\omega$?
    \item [$(\gimel)$] What is the consistency strength of the negative relation at every strong limit singular cardinal $\lambda$, in a universe in which $2^\lambda>\lambda^+$ at every such $\lambda$?
  \end{enumerate}
\end{question}

\newpage

\section{A negative relation from stick}

In the previous sections we used pcf assumptions in order to cope with the problem of Erd\H{o}s and Hajnal.
We move now to the second approach in which prediction principles play an important role.
The prediction principle that we need for the combinatorial proof is called \emph{stick}.
The idea of stick as a prediction principle is well-articulated in \cite[Chapter 4(12)]{hoshea}: ``It consults its stick, its rod directs it''.
Here we need the mathematical incarnation of this idea.
We commence with the following definition.

\begin{definition}
  \label{defstick} Suppose that $\kappa\leq\lambda$.
  \begin{enumerate}
    \item [$(\aleph)$] $\stick(\kappa,\lambda)=\min\{|\mathcal{F}|\mid \mathcal{F}\subseteq[\lambda]^\kappa\wedge \forall y\in[\lambda]^\lambda\exists x\in\mathcal{F}, x\subseteq{y}\}$.
    \item [$(\beth)$] Denote $\stick(\lambda,\lambda^+)$ by $\stick(\lambda)$.
  \end{enumerate}
\end{definition}

The stick principle is closely related to the club principle, but no stationary sets are involved in the prediction.
This fact makes $\stick(\lambda)$ very useful when $\lambda$ is a singular cardinal.\footnote{Notice that $\stick(\lambda)$ denotes both the cardinal value and the combinatorial principle. The principle says that there exists a family of sets $\mathcal{F}\subseteq[\lambda^+]^\lambda$ so that for every $y\in[\lambda^+]^{\lambda^+}$ there exists $x\in\mathcal{F}$ satisfying $x\subseteq{y}$. We adhere to the common usage of this symbol in the literature.}
Let us recall the definition of the club principle (or \emph{tiltan}), which appeared for the first time in \cite{MR438292}.
If $\kappa=\cf(\kappa)>\aleph_0$ and $S\subseteq{\kappa}$ is stationary, then a tiltan sequence $(T_\delta\mid\delta\in{S})$ is a sequence of sets, where $T_\delta$ is a cofinal subset of $\delta$ for each $\delta\in{S}$,\footnote{We assume, tacitly, that $S$ consists of limit ordinals. There is no loss of generality here since $S$ is stationary.} and if $A\in[\kappa]^{\kappa}$ then $S_A=\{\delta\in{S}\mid T_\delta\subseteq{A}\}$ is stationary.
One says that $\clubsuit_S$ holds if there exists such a sequence.

We proceed to the combinatorial result.
Our goal is to prove the relation $\lambda^+\nrightarrow(\lambda^+,(3)_{\cf(\lambda)})^2$ from the stick principle $\stick(\lambda)$.
Negative partition relations at successors of a regular cardinal $\kappa$ follow from $\stick(\kappa)$ as shown in \cite{MR4068775}.
Here we apply a similar idea to successors of singular cardinals.
We need the following lemma about free sets from \cite{MR0424569}.
The lemma and its proof also appear in \cite[Lemma 20.3]{MR795592}.

\begin{lemma}
  \label{lemhm75} Let $\kappa$ be a regular cardinal.
  Suppose that $E=\bigcup_{\alpha\in\kappa}E_\alpha$, and $|E_\alpha|>\kappa$ for every $\alpha\in\kappa$.
  Assume further that $f:E\rightarrow\mathcal{P}(E)$ is a set mapping, and $|f(x)\cap E_\alpha|<\kappa$ for every $x\in{E},\alpha\in\kappa$.
  Then there exists a free set $X$ for $f$ so that $X\cap{E_\alpha}\neq\varnothing$ for every $\alpha\in\kappa$.
\end{lemma}

We can state now the following:

\begin{theorem}
  \label{thmmt} Suppose that $\theta=\cf(\lambda)<\lambda$ and assume that $\stick(\lambda)$ holds.
  Then $\lambda^+\nrightarrow(\lambda^+,(3)_{\cf(\lambda)})^2$.
\end{theorem}

\par\noindent\emph{Proof}. \newline
Let $(\kappa_i\mid 1\leq{i}\in\theta)$ be an increasing sequence of infinite cardinals such that $\lambda=\bigcup_{i\in\theta}\kappa_i$.
Notice that the enumeration of these cardinals begins with $\kappa_1$ since we wish to save the first color to the full-sized independent subsets of the graph.
We shall define a partition $\mathcal{P}=(\mathcal{P}_i\mid i\in\theta)$ of $[\lambda^+]^2$.\footnote{The elements of this partition are not required to be disjoint, so we use here the term \emph{partition} in an unusual way.} Then, essentially, for $\alpha<\beta<\lambda^+$ we will set $c(\alpha,\beta)=i$ iff $\{\alpha,\beta\}\in\mathcal{P}_i$.\footnote{We say \emph{essentially} since the elements of the partition here are not necessarily disjoint, so the formal definition of the coloring will take the first $i$ for which $\{\alpha,\beta\}\in\mathcal{P}_i$.}
The partition $\mathcal{P}$ will be based on a sequence of set-mappings in the following way.
For every $i\in(0,\theta)$ we shall define $f_i:\lambda^+\rightarrow[\lambda^+]^{\leq\kappa_i}$ such that $f_i(\alpha)\subseteq\alpha$ for every $i\in\theta,\alpha\in\lambda^+$.
We let $\mathcal{P}_i=\{\{\alpha,\beta\}\mid\alpha<\beta<\lambda^+,\alpha\in f_i(\beta)\}$.
This procedure defines $\mathcal{P}_i$ for $i>0$, and we let $\mathcal{P}_0=[\lambda^+]^2-\bigcup\{\mathcal{P}_i\mid 1\leq{i}<\theta\}$.

The construction of each $f_i$ is by induction on $\alpha\in\lambda^+$, where at the $\alpha$th stage, $f_i(\alpha)$ is defined simultaneously for each $i\in(0,\theta)$.
Fix $\alpha\in\lambda^+$ and suppose that $f_i(\gamma)$ is already defined for every $\gamma\in\alpha$ and every $i\in\theta$.
Let $(T_\eta\mid\eta\in\lambda^+)$ be a $\stick(\lambda)$ sequence, so $T_\eta\in[\lambda^+]^\lambda$ for every $\eta\in\lambda^+$.
Let $\mathcal{S}_\alpha=\{T_\eta\mid\eta\in\alpha,T_\eta\subseteq\alpha\}$.
Notice that $|\mathcal{S}_\alpha|\leq|\alpha|\leq\lambda$ and hence there exists a decomposition of the form $\mathcal{S}_\alpha=\bigcup\{\mathcal{S}^\alpha_i\mid 1\leq{i}\in\theta\}$, where $i<j\Rightarrow\mathcal{S}^\alpha_i\cap\mathcal{S}^\alpha_j=\varnothing$ and $|\mathcal{S}^\alpha_i|\leq\kappa_i$ for every $i\in(0,\theta)$.

In order to define $f_i(\alpha)$ for each $i\in(0,\theta)$, fix an ordinal $i$ and apply Lemma \ref{lemhm75}, where $\kappa_i^+$ here stands for $\kappa$ there, and $f_i\upharpoonright\alpha$ here stands for $f$ there.
Notice that $|f_i(\gamma)|\leq\kappa_i$ for each $\gamma\in\alpha$ by the induction hypothesis, so the assumptions of the lemma hold.
By the conclusion of the lemma, there exists a free set $X=X_{\alpha{i}}$ for $f_i\upharpoonright\alpha$ which satisfies $X\cap{T}\neq\varnothing$ for every $T\in\mathcal{S}^\alpha_i$.
By removing elements from $X$ if needed, we may assume that $|X|\leq|\mathcal{S}^\alpha_i|\leq\kappa_i$, so we can define $f_i(\alpha)=X=X_{\alpha {i}}$.
This completes the definition of our set mappings, and consequently the definition of $\mathcal{P}$, the partition of $[\lambda^+]^2$.

We define, at this stage, the coloring $c:[\lambda^+]^2\rightarrow\theta$ by letting $c(\alpha,\beta)=i$ iff $i\in\theta$ is the first ordinal so that $\{\alpha,\beta\}\in\mathcal{P}_i$.
We claim that $c$ witnesses the negative relation to be proved.
To see this, let us show firstly that there are no $\alpha<\beta<\delta<\lambda^+$ and $i\in\theta$ such that $c(\alpha,\beta)=c(\alpha,\delta)=c(\beta,\delta)=i$ where $i>0$.
Indeed, if $\alpha<\beta<\delta<\lambda^+$ and $c(\alpha,\delta)=c(\beta,\delta)=i$ then $\{\alpha,\delta\},\{\beta,\delta\}\in\mathcal{P}_i$.
This means that $\alpha,\beta\in f_i(\delta)=X$.
But $X$ is a free set with respect to $f_i\upharpoonright\delta$, and $\beta\in{X}$, hence $f_i(\beta)\cap{X}=\varnothing$.
Since $\alpha\in{X}$ one concludes that $\alpha\notin f_i(\beta)$.
Therefore, $c(\alpha,\beta)\neq{i}$.

Secondly, we argue that there is no $0$-monochromatic subset of $\lambda^+$ of size $\lambda^+$.
To see this, fix $A\in[\lambda^+]^{\lambda^+}$.
Choose an ordinal $\eta\in\lambda^+$ such that $T_\eta\subseteq{A}$.
If $\xi>\eta$ and $\xi>\sup(T_\eta)$ then, by definition, $T_\eta\in\mathcal{S}_\xi$.
Since $A$ is unbounded in $\lambda^+$, one can choose $\xi>\eta,\sup(T_\eta)$ such that $\xi\in{A}$.
Recall that we had a partition $\mathcal{S}_\xi=\bigcup\{\mathcal{S}^\xi_i\mid 1\leq{i}\in\theta\}$, hence $T_\eta\in\mathcal{S}^\xi_i$ for some $i\in(0,\theta)$.

By the choice of $f_i(\xi)$ we know that $T_\eta\cap f_i(\xi)\neq\varnothing$, so one can choose $\alpha\in T_\eta\cap f_i(\xi)$.
The fact that $\alpha\in f_i(\xi)$ implies that $\{\alpha,\xi\}\in\mathcal{P}_i$.
Hence $c(\alpha,\xi)\neq{0}$.
Since $\alpha\in{T_\eta}\subseteq{A}$ and $\xi\in{A}$, one concludes that $c''[A]^2\neq\{0\}$, so we are done.

\hfill \qedref{thmmt}

We would like to emphasize an important aspect of the proof.
The size of each $T_\eta$ is $\lambda$, as our guessing principle is $\stick(\lambda)=\stick(\lambda,\lambda^+)$.
Life would be much simpler if we could replace $\stick(\lambda)$ by $\stick(\kappa,\lambda^+)$ or by $\clubsuit_{S^{\lambda^+}_\kappa}$ for some $\kappa<\lambda$.
One should ask, therefore, whether $\stick(\lambda)$ is essential for the combinatorial proof above.
To wit, one should ask why do we insist on guessing sets of size $\lambda$.
The point lies in Lemma \ref{lemhm75}.
In order to find a free set $X$ that meets every $E_\alpha$, one must verify that $E_\alpha$ is sufficiently large.
In the context of our proof, one needs $|T_\eta|>\kappa_i$ for some relevant $\kappa_i<\lambda$, where $T_\eta$ plays the role of $E_\alpha$ within the proof of the theorem.
However, we do not know in advance the identity of $\kappa_i$, and it might appear as any $\kappa_i$ in the sequence $(\kappa_i\mid i\in\theta)$.
It seems that the only way to cope with this problem is by taking the $T_\eta$'s to be of cardinality $\lambda$.

Though we do not know how to force $\stick(\lambda)$ with the failure of \textsf{SCH} at $\lambda$, we can describe a possible strategy towards this goal.
Let $\lambda$ be a supercompact cardinal and let $\mathscr{U}$ be a normal ultrafilter over $\lambda$.
One says that ${\rm Gal}(\mathscr{U},\lambda^+,\lambda^+)$ holds if every family $\{A_\alpha\mid\alpha\in\lambda\}\subseteq\mathscr{U}$ contains a subfamily $\{A_{\alpha_i}\mid i\in\lambda^+\}$ so that $\bigcap_{i\in\lambda^+}A_{\alpha_i}\in\mathscr{U}$.
One can force such an ultrafilter over a supercompact cardinal $\lambda$, see e.g. \cite{MR4611828}.
If one forces Prikry through such an ultrafilter then every new set of size $\lambda^+$ in the Prikry generic extension contains an old set of size $\lambda^+$.

Indeed, let $\name{A}$ be a name of a set of cardinality $\lambda^+$ and let $p\in\mathbb{P}$ be an arbitrary condition that forces this fact.
For every $\alpha\in\name{A}[G]$ let $p_\alpha=(s_\alpha,A_\alpha)$ be such that $p_\alpha\Vdash\check{\alpha}\in\name{A}$.
Let $A_0\in[A]^{\lambda^+}$ be such that $s_\alpha=s$ for every $\alpha\in{A_0}$.
By ${\rm Gal}(\mathscr{U},\lambda^+,\lambda^+)$ there exists $B\in\mathscr{U}$ and there is $A_1\in[A_0]^{\lambda^+}$ such that $B\subseteq{A_\alpha}$ for every $\alpha\in{A_1}$.
Set $q=(s,B)$ and notice that $p_\alpha\leq{q}$ for every $\alpha\in{A_1}$.
Now $A_1\in{V}$ and $q\Vdash A_1\subseteq\name{A}$.
For a comprehensive account of density properties with respect to Prikry forcing, see \cite{MR3577887}.

Thus if stick or tiltan hold at $\lambda^+$ in the ground model, this will be preserved in the generic extension.
Moreover, this is true for a variety of Prikry-type forcing notions, including Prikry forcing with interleaved collapses.
There are also several ways to force stick or tiltan at a supercompact cardinal $\lambda$ while increasing $2^\lambda$ above $\lambda^+$.
We do not know, however, how to force these two things together:

\begin{question}
  \label{qstick} Is it consistent that $\lambda$ is supercompact, $\mathscr{U}$ is a normal ultrafilter over $\lambda$ satisfying ${\rm Gal}(\mathscr{U},\lambda^+,\lambda^+)$, $\stick(\lambda)$ holds and $2^\lambda>\lambda^+$?
\end{question}

\newpage

\section{Appendix}

In this section we describe the extender-based Prikry forcing with interleaved collapses, and we prove some basic facts with respect to this forcing notion.
The main statement to be proved is that this forcing notion enjoys the strong Prikry property.
By and large, the results presented in this section are known.
However, some of the statements and proofs are folklore, and apparently there is no citable source for them.
The presentation of the extender-based Prikry forcing in this section follows the conventions set by Merimovich in various articles.

\begin{definition}[Strong cardinal]
\label{defstrongcardinal}
Let $\kappa$ be a cardinal. We say that $\kappa$ is $\beta$-strong if there is an elementary embedding $j\colon V \to M$, $M$ is closed under $\kappa$-sequences, $\crit j = \kappa$ and $V_{\beta}\subseteq M$. We shall say that $\kappa$ is strong if $\kappa$ is $\beta$-strong for every $\beta$.
\end{definition}

Thus, the concept of a strong cardinal depicts a global property.
This property comes in handy while trying to use extenders.
Fix an elementary embedding $j\colon V \to M$ with critical point $\kappa$ and $\alpha \leq j(\kappa)$. Following Merimovich, let us derive an extender $E$ by assigning for each $d \in [\alpha\setminus \kappa]^{\leq \kappa}$, the ultrafilter $E_d = \{X \subseteq (d\times \kappa)^{<\kappa} \mid (j\restriction d)^{-1} \in j(X)\}$.
Note that for $d = \emptyset$, $E_d$ is principal so we typically remove it from the extender (even though this will not play any role in the construction).

Let $d \subseteq e$ be sets in $[\alpha \setminus \kappa]^{\leq\kappa}$. Then, there is a Rudin-Keisler projection $\pi_{e,d}$ between $E_e$ and $E_d$ given by $x\mapsto x\restriction d$.
We identify $E$ with the pair $\langle \langle E_d \mid d \in [\alpha\setminus \kappa]^{\leq \kappa}\rangle, \langle \pi_{e,d} \mid d \subseteq e \subseteq \alpha\setminus \kappa,\, |e| \leq \kappa\rangle\rangle$ and say that $E$ is a $(\kappa,\alpha)$-extender derived from $j$. In particular, if $\kappa$ is $\kappa+2$-strong then there is a $(\kappa,\kappa^{++})$-extender.

Given an extender one can take the ultrapower by the extender in order to derive an elementary embedding, from which the extender can be recovered. Note that the assumption that the extender is derived from an elementary embedding with critical point $\kappa$ entails some combinatorial properties that we did not explicitly state in the definition of the extender.

\begin{lemma}
\label{lembasic1}
Let $E$ be as above.
\begin{enumerate}
\item [$(\aleph)$] For every $d$, the map $x \mapsto x(\kappa)$ from $(d\times \kappa)^{<\kappa}$ is well defined on an $E_d$-large set and is equal to its Rudin-Keisler projection to the normal measure.
\item [$(\beth)$] Let $d \in [\alpha\setminus\kappa]^{\kappa}$ and let $g\colon \kappa \to d$ be a bijection. Then, there is a large set $B\in E_d$ such that for every $x\in B$, $g\restriction x(\kappa)$ is a bijection from $x(\kappa)$ to $\dom x$.
\end{enumerate}
\end{lemma}

\par\noindent\emph{Proof}. \newline
Both parts follow directly from the definition. The second part plays an important role in \cite{MR4264147}.

\hfill \qedref{lembasic1}

While the ultrafilters $E_d$ are not normal, they support a type of diagonal intersection.

\begin{lemma}
\label{lemdiagonal}
Let $A_* = \langle A_\eta \mid \eta \in B\rangle$ where $B\in E_d$. We assume that $B \subseteq [d\times \kappa]^{<\kappa}$, consists of increasing sequences and respects some fixed enumeration of $d$. Let
\[\triangle^\pi A_\eta = \{\zeta\in B \mid \forall \eta \in B, \sup \range \eta < \zeta(\kappa) \Rightarrow \zeta \in A_\eta\}\]
Then, $\triangle^\pi A_\eta\in E_d$.
\end{lemma}

\par\noindent\emph{Proof}. \newline
Let us verify that the seed of $E_d$, $j_E^{-1}\restriction (j\image d)$ belongs to $j_E(A_*)$.
Note that \[\{\eta \in j_E(B) \mid \sup\range \eta < j_E^{-1}\restriction (j\image d)(\kappa)=\kappa\} = j\image B.\]
Indeed, every such element must be smaller than $\kappa$ and consequently its domain is contained in $j(g)\image \kappa = j\image d$.
Thus, the lemma follows from the definition of $E_d$.

\hfill \qedref{lemdiagonal}

Given an extender $E$, we define the Extender based Prikry forcing, following \cite{MR4264147}. We would like to incorporate collapses as well, and in the general case there are some subtleties, which are pcf-related, that limit the general possibilities.

Let us assume from now on that $E$ is a $(\kappa, \kappa^{++})$-extender.
Let $h \colon \kappa \to \kappa$ be such that for all $\alpha$, $h(\alpha) > \alpha^{++}$ and $h(\alpha)$ is a regular cardinal. We would like to define the extender based Prikry forcing with interleaved collapses, based on $h$. The following statements should hold in the generic extension:

\begin{itemize}
\item $\cf \kappa = \omega$ as witnessed by a Prikry sequence $\langle \rho_n \mid n < \omega\rangle$,
\item $2^\kappa \geq \kappa^{++}$,
\item Every cardinal in the open intervals $(h(\rho_n), \rho_{n+1})$ as well as $(h(\omega), \rho_0)$ is collapsed and every other cardinal is preserved.
\item If we start with a model of $\GCH$, then $\GCH$ holds below $\kappa$ in the generic extension.
\end{itemize}

\begin{lemma}
\label{lemgenericcummings}
Assume $2^\kappa = \kappa^{+}$. Let $j\colon V \to M$ be the ultrapower by the extender $E$. Let $\iota \colon V \to N$ be the ultrapower by the normal measure of $j$ and let $k$ be the factor map.
Then, there is an $N$-generic filter $\tilde K\subseteq {\rm Col}^{N}(\iota(h)(\kappa), <\iota(\kappa))$ and moreover $K = k\image \tilde K$ generates an $M$-generic filter.
\end{lemma}

\par\noindent\emph{Proof}. \newline
See \cite[Lemma 15.3]{MR2768691}.
\hfill \qedref{lemgenericcummings}

\begin{lemma}
\label{lemmutualgeneric}
Let $K$ be as above. Then, the filters $\langle K, j_E(K), \dots, j_E^m(K)\rangle$ are mutually generic over the $m+1$-iteration of the ultrapower by $E$.
\end{lemma}

\par\noindent\emph{Proof}. \newline
First notice that the filters $\tilde{K}, \iota(\tilde{K}), \dots, \iota^m(\tilde{K})$ are mutually generic by the Cummings-Woodin lemma from \cite[Fact 2]{MR1041044}. Second, by induction, let us lift those filters to the iterated ultrapower by $E$. This can be done since the width of the maps is below the distributivity of the forcing notions.

\hfill \qedref{lemmutualgeneric}

Note that the filter $\langle K, j_E(K), \dots, j_E^m(K)\rangle$ is generated by elements of the form $\langle j_E(C)(\kappa), j^2_E(C)(j_E(\kappa)), \dots, j^m_E(C)(j^{m-1}_E(\kappa))\rangle$.
Let us define the forcing notion. Following Merimovich, a condition in the forcing is
\[p = \langle f, A, c_0,\dots, c_n, C\rangle\]
where:
\begin{enumerate}
\item $f$ is a partial function from $\kappa^{++} \setminus \kappa$ to $\kappa^{<\omega}$ with $|\dom f| < \kappa$ and $\kappa \in \dom f$.\footnote{We will identify the functions $f$ with conditions in $\Add(\kappa,\kappa^{++})$.}
\item Let $d = \dom f$, then $A \in E_d$. For every $h \in A$ we assume that $h$ is strictly increasing, $h(\kappa)$ is inaccessible and larger than every ordinal that appears in $f(\alpha)$ for $\alpha \in d$.
\item $n = \dom f(\kappa)$ and let $\langle \rho_i \mid i < n\rangle = f(\kappa)$. We denote $\rho_{-1} = \omega$ and $\rho_{n} = \kappa$.
\item For every $i \leq n$, $c_i \in {\rm Col}(h(\rho_{i-1}), <\rho_i)$.
\item $C$ is a function with domain $A$ with the property that for every $x, y\in A$ if $x(\kappa) = y(\kappa)$ then $C(x) = C(y)$. Moreover, $[C, d]_{E} \in K$.
\end{enumerate}

Given conditions $p = \langle f^p, A^p, c_0^p,\dots, c_n^p, C^p\rangle$, $q = \langle f^q, A^q, c_0^q,\dots, c_m^q, C^q\rangle$ we say that $p \leq^* q$ (i.e.\ $q$ is a direct extension of $p$) if:
\begin{enumerate}
\item $n = m$,
\item $f^q \supseteq f^p$.
\item $\pi_{d^q, d^p}(A^q) \subseteq A^p$.
\item For every $i \leq n$, $c_i^q \supseteq c_i^p$.
\item For every $x \in A_q$, $C^q(x) \supseteq C^p(\pi_{d^q,d^p}(x))$.
\end{enumerate}

Let $p = \langle f^p, A^p, c_0^p,\dots, c_n^p, C^p\rangle$ be a condition and $\eta \in A^p$.
The condition $q = p^\smallfrown \eta$ is defined by
\begin{itemize}
\item $\dom f^q = \dom f^p$. For every $\alpha \in \dom \eta$, $f^q(\alpha) = f^p(\alpha)^\smallfrown \langle \eta(\alpha)\rangle$, and otherwise $f^q(\alpha)=f^p(\alpha)$.

In particular, the length of $q$ is $n+1$.
\item $A^q = \{x \in A^p \mid x(\kappa) > \sup\range \eta\}$.
\item For every $i \leq n$, $c_i^p = c_i^q$.
\item $c_{n+1}^q = C(\eta)$.
\item $C^q = C^p \restriction A^q$.
\end{itemize}
We say that $q$ is a one-step extension of $p$.

Let us define recursively that $q$ is an $n$-step extension by setting $p$ as the only $0$-step extension of $p$ and letting $q$ be an $n$-step extension of $p$ iff it is a one-step extension of an $(n-1)$-step extension of $p$.
Finally, we define $p \leq q$ if there is a finite sequence of direct extensions and one-step extensions leading from $p$ to $q$.

\begin{lemma}
\label{lemtransitive}
$\leq^*$ is a transitive relation.
\end{lemma}

\begin{lemma}
\label{lempurextension}
$p\leq q$ if and only if there is $r$ which is a $k$-step extension of $p$ such that $r\leq^* q$.
\end{lemma}

\par\noindent\emph{Proof}. \newline
Let us prove that given conditions $p,q$ and $r$ such that $p \leq^* r$ and $q = r^\smallfrown \eta$ is a one-step extension of $r$, one can find $r'$ such that $r'$ is a one-step extension of $r$ and $r'\leq^* q$.

Let $n$ be the length of $f^q(\kappa)$.
Let $f^{r'} = f^r\restriction (d - \dom \eta) \cup f^q\restriction (\dom \eta)$.
Let $A^{r'} = A^r \cup \{\eta\}$, $c_i^{r'} = c_i^r$ for $i \leq n$, and $C^{r'} = C^r \cup \{\langle \eta, c_{n+1}^r\rangle\}$.
Note that for every $x\in A^r$, $x(\kappa)$ is larger than $\eta(\kappa)$ and thus requirement (5) holds.
Repeating this process for finitely many times we obtain the desired conclusion.

\hfill \qedref{lempurextension}

This forcing is a typical example of a Prikry type forcing: it carries two types of partial orders, the direct order which satisfies better closure properties and the actual order with which we will force. In order to gain better control over the properties of the extension we would like to decompose the direct order $\leq^*$ into components.

We say that $p\leq^{**} q$ if $p\leq^* q$ and $c_i^p = c_i^q$ for all $i$.

We say that $p\leq^{col} q$ if $p \leq^* q$ and $f^p = f^q, A^p = A^q, C^p = C^q$.

Our next goal is to prove the strong Prikry property of the above defined forcing notion.
In fact, we shall prove the \emph{complete Prikry property}, defined in \cite{MR4325864}.
In that paper, it was illustrated that this property gives a cleaner way to understand the properties of a Prikry type forcing.
It was also demonstrated in general that the complete Prikry property implies the strong Prikry property.

\begin{lemma}
\label{lemspp}
Let $D\subseteq \mathbb{P}$ be dense open and let $p$ be a condition. Then, there is a condition $q$ and a natural number $n$ such that:
\begin{itemize}
\item $p \leq^{**} q$.
\item There is a dense open set in $\leq^{col}$ above $q$, $D_{col}$, such that for every $r \in D_{col}$, every $n$-step extension of $r$ belongs to $D$.
\end{itemize}

Moreover, the \emph{Complete Prikry Property} holds: for every $\leq^*$-open set $U$, there is a condition $q$ such that $p\leq^{**} q$ and for every natural number $n$:
\begin{itemize}
\item either for every $n$-step extension $r$ of $q$ there is no $\leq^{col}$ extension of $r$ in $U$ or
\item there is a dense subset in $\leq^{col}$ above $q$, such that for every $q'\in D_{col}$ and for every $n$-step extension of $q'$,  $r$, we have $r \in U$.
\end{itemize}
\end{lemma}

\par\noindent\emph{Proof}. \newline
Let $U$ be $\leq^*$-open.
Let $\chi$ be a regular cardinal larger than $(2^\kappa)^{++}$.
Fix a condition $p$ and fix an elementary substructure $M \prec H(\chi)$ containing $p, U,\mathbb{P}$, of cardinality $\kappa$, which is closed under $<\kappa$-sequences.

Let us enumerate all dense open sets in the forcing $\Add(\kappa,\kappa^{++})$ in $M$ by $\langle \tilde{D}_{\alpha} \mid \alpha < \kappa\rangle$. Let us define a sequence of Cohen conditions $f_\alpha$ as follows.
Let $f_0 = f^p$ and for $0 < \alpha < \kappa$ limit, let $f_\alpha = \bigcup_{\beta < \alpha} f_\beta$.
For the successor stage, let us assume that $f_\alpha$ was defined and let $d_\alpha = \dom f_\alpha$.
Let $f_{\alpha+1} \in \tilde{D}_{\alpha}\cap M$ be stronger than $f_\alpha$.
Finally, let $f_\kappa = \bigcup_{\alpha < \kappa} f_\alpha$, and $d_\kappa = \dom f_\kappa$. Let us define for every $\alpha \leq \kappa$,
\[p_\alpha= \langle f_\alpha, \pi^{-1}_{d_\alpha, d_0}(A^p), c^p_0,\dots, c^p_n, C^p\circ \pi_{d_\alpha, d_0}\rangle,\]
let $A^\alpha = \pi^{-1}_{d_\alpha, d_0}(A^p)$.

Let $\vec\eta = \langle \eta_0,\dots, \eta_{m-1}\rangle$ be in $\pi^{-1}_{d_\alpha, d_0}(A^p)$ such that the $m$-step extension of $p_\alpha$ by $\vec\eta$ is a condition.
Let us look at all extensions $g$ of $f_\alpha$ such that either there is a direct extension $r$ of $p_\alpha^\smallfrown \vec\eta$ in $U$ such that $f^r = g^\smallfrown \vec\eta$ or that for any further extension of $g$ there is no such direct extension $r \in U$.
This is clearly a dense open set in $M$ and thus $f_\kappa$ belongs to any such set.
In particular, for every $\vec\eta$ such that there is a direct extension of $p_\kappa^\smallfrown \vec\eta$ in $U$, such a direct extension will not modify the Cohen part $f^{p_\kappa}$.

Let us now deal with the collapse parts.
Let $\mathcal{D}$ be the collection of all $[C]_E \leq [C^\kappa]_E$ such that for measure one many $\eta \in A^\kappa$, there is a direct extension $p_\kappa^\smallfrown \eta \leq^* r$ such that $c_{n+1}^r = C(\eta)$ and no further direct extension $r \leq^* s$ belongs to $U$.
Note that $\mathcal{D}$ is open. 
If there is $[C]_E\in K\cap \mathcal{D}$ above $[C_\kappa]_{E}$ (e.g., if $\mathcal{D}$ is dense), let us extend $p^\kappa$ by replacing $C_\kappa$ with $C$ and shrink $A^\kappa$ to the large set.

Now, for every $\eta \in A^\kappa$ there is an extension of the first $n$ components of the collapse part of $p_\kappa^\smallfrown \eta$ such that no further direct extension enters $U$. By shrinking $A^\kappa$ again\footnote{The lower collapses corresponding to $\eta$ are in $V_{\eta(\kappa)}$ and in the ultrapower then represent a member of $V_\kappa$, which is the corresponding value obtained on a large set.}, we can stabilize this extension and by strengthening $p_\kappa$ we may assume that it was the collapse part to begin with.

Let us assume that for every step up to $m$ we constructed the relevant $C$. At step $m$, let us consider the set $\mathcal{D}_m$ that consists of all functions $\langle C_0, C_1, \dots, C_m\rangle$ with $\dom C_i = [A^\kappa]^{i+1}$ such that on some large subset of $A^\kappa$ for every $\vec\eta = \langle \eta_0,\dots, \eta_m\rangle$ increasing\footnote{Namely $\sup \range \eta_i < \eta_{i+1}(\kappa)$.}, there is a direct extension $p_\kappa^\smallfrown \vec\eta \leq^* r$ with collapse parts $C_0(\eta_0), C_1(\eta_0, \eta_1)\dots, C(\eta_0,\dots, \eta_{n-1})$, for which no further direct extension belongs to $U$.
If there is $[C] \in K$ above $[C^\kappa]$ such that taking $C_i$ as $C(\eta_i(\kappa))$, the obtained condition meets $\mathcal{D}$, we replace $[C^\kappa]$ by this $C$ and shrink the large sets accordingly.
Eventually, there must be $m$ for which we cannot find such $C$.
Otherwise, by the $\sigma$-closure of the direct extension after $\omega$ many steps we will have no extension of the condition in $U$.

Let us assume that at step $m$ there is no such $C$.
So, necessarily, there is $C$ in $K$ such that no strengthening of $C$ belongs to $\mathcal{D}_m$. So, for every strengthening of $C$ there is a further extension of the collapse part that enters $U$.
Thus, we obtain a dense set below $[C^\kappa]$ and therefore there is such a condition in $K$.

Look at all extensions of the first $n$ collapses that enter $U$. Taking for each $\eta$ a maximal antichain we again obtain a regressive function that can be stabilized on a large set. Taking a $\pi$-diagonal intersection of the large sets, we shrink the large set as required.

\hfill \qedref{lemspp}

From the strong Prikry property one can derive the following version of the Prikry property.
\begin{lemma}
\label{lemversionofpp}
Let $\Phi$ be a statement in the language of forcing and $p$ a condition.

Then, there is a condition $q$ such that $p\leq^{**} q$ and there is a dense subset $D_{col}$ in $\leq^{col}$ above $q$ of conditions that decide the truth value of $\Phi$.
\end{lemma}

These properties can be used to compute which cardinals are collapsed in the generic extension.

\begin{lemma}
\label{lemtcfsequences}
Assume that $2^\kappa=\kappa^+$.
Let $G\subseteq \mathbb{P}$ be $V$-generic and let $g\colon \kappa \to \kappa$ be a function in $V$ such that for every inaccessible cardinal $\alpha$, $g(\alpha) \geq \alpha$ is a regular cardinal. Let us assume further that in $V[G]$, $g(\rho_n)$ is regular for every $n$ (so $g(\alpha) \leq h(\alpha)$ for every $\alpha$).
Then, if $g(\alpha) > \alpha^{++}$ on a measure one set, \[\mathrm{tcf} \prod g(\rho_n) = \kappa^{+}.\]
\end{lemma}

\par\noindent\emph{Proof}. \newline
Let $\name{\ell}$ be a name for a member of $\prod g(\name{\rho}_n)$ and let $p$ be a condition.
Using the strong Prikry property, and the $\sigma$-closure of the direct extension, there is a condition $q$ such that for every $n$ there is $m$ and a dense subset of the collapse parts $D_n$ such that every $m$-step extension of a $\leq^{col}$-extension of $q$ from $D_n$ decides the value of $\ell(n)$.

Without loss of generality, $p$ is of length $0$. Using the $\sigma$-closure of the collapses, we may assume that $q$ belongs to $D_n$ for all $n$. As $\Vdash \ell(n) \leq h(\rho_n) < \rho_{n+1}$, we may assume that $m = n$.
Thus, taking $d = \dom q$, $\ell(n)$ is a function of the first $n$ many Prikry points of $E_d$.

Let us now show that there is an $E_d$-measure-one set $B$ such that for every $\rho < \kappa$, $\{y \in B \mid y(\kappa) = \rho\}$ is at most $\rho^{++}$.
Indeed, the seed of $E_d$, $s_d = j^{-1}\restriction j\image d$ satisfies that its range is bounded by $\kappa^{++} = s_d(j(\kappa))^{++}$. As the domain is fixed, the number of possible values is at most $|(\rho^{++})^{\rho}| = \rho^{++}$.
Thus, one can bound $\ell$ based only of the values of $\langle \rho_n \mid n < \omega\rangle$. So, there is a function $f_{\ell}\colon \kappa \to \kappa$ such that $q \Vdash \forall n, \ell(n) \leq f_\ell(\rho_n)$.
By the assumption of $2^\kappa=\kappa^+$ in the ground model, we obtain a family of functions of size $\kappa^+$ that dominates every function in the scale.
Thus, the length of the scale is $\kappa^{+}$.

\hfill \qedref{lemtcfsequences}

\section{Acknowledgements}
The authors would like to thank the anonymous referee for their useful suggestions and corrections that greatly improved this manuscript.

\newpage

\bibliographystyle{alpha}
\bibliography{arlist}

\end{document}